\documentclass[a4paper,10pt,reqno]{amsart}
\usepackage[centertags]{amsmath}
\usepackage{amsfonts,amssymb,amsthm}
\usepackage{newlfont}
\usepackage[dvips]{graphicx} %allows.eps and.epsi graphics to be inserted
\usepackage[english]{babel}
\usepackage[body={15cm, 23cm},centering]{geometry}
\usepackage{latexsym}
\usepackage{amsmath}
\usepackage{amsfonts}
\usepackage{amssymb}
 \theoremstyle{plain}
 \begingroup
 \newtheorem{teo}{Theorem}[section]
 \newtheorem{lem}[teo]{Lemma}
 \newtheorem{prop}[teo]{Proposition}
 \newtheorem{cor}[teo]{Corollary}
 \endgroup

 \theoremstyle{definition}
 \begingroup

 \endgroup

 \theoremstyle{remark}
 \begingroup
 \newtheorem{rem}[teo]{Remark}
 \endgroup

 \numberwithin{equation}{section}

\newcommand{\D}{\Delta}

\def\s{\,\,\,\,}

\def\dint{\displaystyle{\int}}

\def\mv{1.7ex}
\newcommand{\R}{{\mathbb R}}

\begin{document}
\title[M-T Problem ]{Extremal functions for Moser-Trudinger type inequality on compact
closed 4-manifolds}
\author[Yuxiang Li, Cheikh Birahim Ndiaye]{Yuxiang Li,  Cheikh Birahim Ndiaye}
\date{\today}

\address[Yuxiang Li]{ICTP, Mathematics Section, Strada Costiera 11, 34014 Trieste, Italy}
\email[yuxiang li]{liy@ictp.it}

\address[Cheikh Birahim Ndiaye]{S.I.S.S.A., Via Beirut 2-4, 34014, Trieste, Italy}
\email[Cheikh Birahim Ndiaye]{ndiaye@sissa.it}

\maketitle
%%-----------------------------------------------------------------------------------------

\begin{center}
\begin{minipage}{12cm}
\small{ \noindent{\bf Abstract.} Given a compact closed four
dimensional smooth Riemannian manifold, we prove existence of
extremal functions for Moser-Trudinger type inequality. The method
used is Blow-up analysis combined with capacity techniques. }
\end{minipage}
\end{center}

\begin{center}
\begin{minipage}{12cm}
\vspace{10pt} \noindent {\bf Keywords}: Moser-Trudinger inequality,
Blow-up analysis, Capacity, Extremal function, Green function.
\end{minipage}

\end{center}
\vspace{10pt} \noindent

\begin{center}
\begin{minipage}{12cm}
\vspace{6pt} \noindent {\bf 2000 Mathematics Subject
Classification}: 46E35, 26D10
\end{minipage}
\end{center}

\bigskip
\tableofcontents
%%--------------------------------------------------------------------------------------------
\section{Introduction}
It is  well-known that Moser-Trudinger type inequalities are crucial
analytic tools in the study of partial differential equations
arising from geometry and physics.

In fact, much work has been done on such inequalities and their
applications in the last decades, see for example, \cite{adr},
\cite{bw}, \cite{csya}, \cite{cy}, \cite{dm}, \cite{mn}, \cite{p},
 and the references therein.

There are two important objects in the study of Moser-Trudinger type
inequalities: one is to find the best constant and the other is to
determine whether there exist extremal functions.

For the best constant there are the celebrated work of
Moser\cite{mo} and the extension to higher order derivatives by
Adams \cite{adr} on flat spaces. In the context of curved spaces
Fontana has extended the results of Adams, see \cite{fl}.

To mention results about extremal functions, we cite the results of
Carleson and Chang \cite{cc}, Flucher \cite{fm} and Lin \cite{lin} in
the Euclidean case and the results of Li \cite{ly1}, \cite{ly2} in
the curved one. In \cite{ly1} and \cite{ly2} the author have proved
the existence of an extremal function for the classic
Moser-Trudinger inequality on a compact manifold under a constraint
involving only the first derivatives.

In this paper, we will extend the results of Li to a compact closed
four dimensional smooth Riemannian manifold under a constraint
involving the Laplacian. More precisely we prove the following
Theorems
\begin{teo}\label{eq:theorem1}
Let \;$(M,g)$\;be a compact closed smooth\;$4$-dimensional
Riemannian manifold. Then setting
\begin{equation*}
\mathcal{H}_{1}=\{u\in H^2(M):\overline u=0,\;\;\int_{M}|\D_g
u|^2dV_g=1\}
\end{equation*}
we have that
\begin{equation*}
\sup_{u\in\mathcal{H}_{1}}\int_{M}e^{32\pi^2u^2}dV_g
\end{equation*}
is attained.
\end{teo}
\vspace{10pt}

On the\;$4-$dimensional manifold \;$(M,g)$\;, the so-called Paneitz
operator, which is defined in terms of the scalar
curvature\;$R_g$\;and the Ricci tensor \;$Ric_g$ as
\begin{equation*}
P^4_gu=\D_g^2u+div_g(\frac{2}{3}R_gg-2Ric_g)du\;\;\;u\in
C^{\infty}(M),
\end{equation*}
plays an important role in conformal geometry see \cite{csya},
\cite{cy}, \cite{cgy2}, \cite{dm}, \cite{g}, \cite{nd}, \cite{p1}.
In particular, the relation between the Paneitz operator and
the \;$Q$-curvature, which is defined as
\begin{equation}
Q_g=-\frac{1}{12}(\D_{g}R_{g}-R_{g}^{2}+3|Ric_{g}|^{2}),
\end{equation}
is of great interest.
It is well-known that Moser-Trudinger inequalities involving
\;$P^4_g$\; play an important role in the problem of prescribing
constant \;$Q$-curvature see \cite{dm}, \cite{Li-Li-Liu}, \cite{nd}.
Therefore it is worth having an extension of
Theorem\;$\ref{eq:theorem1}$\; concerning the Paneitz operator as
well. Our next result goes in this direction. More precisely we have
the following.
\begin{teo}\label{eq:theorem2}
Let \;$(M,g)$\;be a compact closed smooth\;$4$-dimensional
Riemannian manifold. Assuming that \;$P^4_g$\;is non-negative and
\;$ker P^4_g\simeq\R$, then setting
\begin{equation*}
\mathcal{H}_{2}=\{u\in H^2(M):\overline u=0,\;\;<P^4_gu,u>=1\}
\end{equation*}
we have
\begin{equation*}
\sup_{u\in\mathcal{H}_{2}}\int_{M}e^{32\pi^2u^2}dV_g
\end{equation*}
is attained.
\end{teo}
\begin{rem}
Since the leading term of \;$P^4_g$( for the definition see the
Section 2)\;is \;$\D_g^2$\;then the two Theorems are quite similar.
We point out that the same proof is valid for both except some
trivial adaptations, hence we will give a full proof of
Theorem\;$\ref{eq:theorem1}$\;only and  sketch the proof of
Theorem\;$\ref{eq:theorem2}$ in the last section.
\end{rem}
\begin{rem}
We mention that due to a result by Gursky, see \cite{g} if both the
Yamabe class \;$Y(g)$\;and\;$\int_{M}Q_gdV_g$\;are non-negative,
 then we have that \;$P^4_g$\;is non-negative and \;$ker P^4_g\simeq\R$.
\end{rem}
We are going to describe our approach to prove
Theorem\;$\ref{eq:theorem1}$. We will use Blow-up analysis. First of
all we take  a sequence \;$(\alpha_k)_k$\; such that
\;$\alpha_k\nearrow32\pi^2$, and by using Direct Methods of the
Calculus of variations we can find \;$u_k\in \mathcal{H}_1$\; such
that
\begin{equation*}
\int_Me^{\alpha_k u_k^2}dV_g=\sup_{v\in \mathcal{H}_1
}\int_Me^{\alpha_kv^2}dV_g.
\end{equation*}
see Lemma\;$\ref{eq:const}$. Moreover using the Lagrange multiplier
rule  we have that \;$(u_{k})_k$\;satisfies the equation:
\begin{equation}\label{eq:edpk}
\Delta_g^2u_k=\frac{u_k}{\lambda_k}e^{\alpha_ku_k^2}-\gamma_k
\end{equation}
for some constants $\lambda_k$ and $\gamma_k$.

Now it is easy to see that if there exists \;$\alpha>32\pi^2$\;such
that $\int_{M}e^{\alpha u_k^2}dV_g$ is bounded, then by using
Lagrange formula, Young's inequality and Rellich compactness
Theorem, we obtain that the weak limit of\;$u_k$\;becomes an
extremizer. On the other hand if
\begin{equation*}
c_k=\max_{M}|u_k|=|u_k|(x_k)
\end{equation*}
is bounded, then from standard elliptic regularity theory \;$u_k$\;is
compact, thus  converges uniformly to an extremizer.
Hence assuming that Theorem\;$\ref{eq:theorem1}$\; does not hold, we get\\
1)
\begin{equation*}
\forall
\alpha>32\pi^2\;\;\;\lim_{k\rightarrow+\infty}\int_{M}e^{\alpha
u_k^2}dV_g\rightarrow+\infty
\end{equation*}
2)
\begin{equation*}
c_k\rightarrow +\infty
\end{equation*}
We will follow the same method as in \cite{ly1}  up to some extents.\\

In \cite{ly1}, the function sequence we studied is the following:
\begin{equation*}
-\Delta_gu_k=\frac{u_k}{\lambda_k}e^{\alpha_k'u_k^2}-\gamma_k,
\end{equation*}
where $\alpha_k'\nearrow 4\pi$, and $u_k$ attains
$\sup\limits_{\int_M|\nabla_gu|^2dV_g=1, \bar{u}=0}
\int_Me^{\alpha_k'u^2}dV_g$. We also assumed $c_k\rightarrow
+\infty$. Then  we have
\begin{equation}\label{l1}
2\alpha_kc_k(u_k(x_k+r_kx)-c_k)\rightarrow -2\log(1+\pi|x|^2)
\end{equation}
for suitable choices of $r_k$, $x_k$. Next we proved the following
\begin{equation}\label{l2}
\lim_{k\rightarrow+\infty}\int_{\{u_k\leq\frac{c_k}{A}\}}|\nabla_gu_k|^2dV_g
=\frac{1}{A}\;\;\forall A>1,
\end{equation}
which implies that
$$\lim_{k\rightarrow+\infty}\dint_Me^{\alpha_ku_k^2}dV_g=\mu(M)+\lim_{k\rightarrow+\infty}
\frac{\lambda_k}{c_k^2},$$ and that $c_ku_k$ converges to some Green
function weakly. In the end, we got an upper bound of
$\frac{\lambda_k}{c_k^2}$ via capacity.

\begin{rem} (\ref{l1}) was first discovered by Struwe in \cite{St}.
\end{rem}

\begin{rem} (\ref{l2}) also appeared in \cite{ao}.\\
\end{rem}

However there are two main differences between the present case and
the one in \cite{ly1}. One is that there is no direct maximum
principle for equation ($\ref{eq:edpk}$) and the other one is that
truncations are not allowed in the space \;$H^2(M)$\;. Hence to get a
counterpart of (\ref{l1}) and ($\ref{l2}$) is not easy.

To solve the first difficulty, we replace $c_k(u_k(x_k+r_kx)-c_k)$
with $\beta_k(u_k(exp_{x_k}(r_kx))-c_k)$, where
$$1/\beta_k=\int_M\frac{|u_k|}{\lambda_k}e^{\alpha_ku_k^2}dV_g.$$
By using the strength of the Green representation formula, we get
that the profile of \;$u_k$\;is either a constant function or a
standard bubble. The second difficulty will be solved by applying capacity
and Pohozaev type identity. In more detail we will prove that
$\beta_ku_k\rightharpoonup G$\;(see Lemma\;$\ref{eq:profile}$) which
satisfies
\begin{equation*}
\left\{\begin{array}{l}
           \Delta_g^2G=\tau(\delta_{x_0}-{Vol_g(M)})\\
           \int_MG=0.
          \end{array}\right.
      \end{equation*}
for some $\tau\in(0,1]$. Then we can derive from  a Pohozaev type
identity (see Lemma\;$\ref{eq:P}$) that
\begin{equation*}
\lim_{k\rightarrow+\infty}\int_Me^{\alpha_ku_k^2}dV_gVol_g(M)+\lim_{k\rightarrow+\infty}\tau^2\frac{\lambda_k}{\beta_k^2}.
\end{equation*}
In order to apply the capacity, we will follow some ideas in \cite{Li-Li-Liu}.
Concretely, we will show that up to a small term  the energy of
\;$u_k$\;on some annulus is bounded below by the Euclidean one (see
Lemma\;$\ref{eq:capeu}$). Moreover one can prove the existence of
\;$U_k$\;(see Lemma\;$\ref{eq:}$)\; such that the energy of
\;$U_k$\;is comparable to  the Euclidean energy of \;$u_k$\;, and
the Dirichlet datum and Neumann datum of \;$U_k$\; at the
boundary of the annulus are also comparable to those of $u_k$.
In this sense, we simplify the calculation of capacity in \cite{ly2}.
Now using
capacity techniques  we get
\;$\frac{c_k}{\beta_k}\rightarrow d$\;and\; $d\tau=1$, see
Proposition\;$\ref{eq:upperbound}$. Furthermore we have that
\begin{equation*}
\lim_{k\rightarrow+\infty}
\tau^2\frac{\lambda_k}{\beta_k^2}\leq\frac{\pi^2}{6}
e^{\frac{5}{3}+32\pi^2S_0}.
\end{equation*}
Hence we arrive to
\begin{equation}
\label{eq:contrad}\sup_{u\in\mathcal
{H}_1}\int_{M}e^{32\pi^2u^2}dV_g\leq
Vol_g(M)+\frac{\pi^2}{6}e^{\frac{5}{3} +32\pi^2S_0}.
\end{equation}
In the end, we will  find test functions  in order to
contradict\;$\eqref{eq:contrad}$. We will simplify the arguments in
\cite{ly1}. Indeed we use carefully the regular part of \;$G$\;to
avoid cut-off functions
and hence making the calculations simpler.\\

The plan of the paper is the following: In Section 2 we collect some
preliminary results regarding the existence of the Green functions
for \;$\D^2_g$\; and \;$P^4_g$\;, and associated Moser-Trudinger
type inequality. In Section 2 we prove Theorem\;$\ref{eq:theorem1}$.
This Section is divided into six subsections. In the first one, we
deal with concentration behavior and the profile of the blowing-up
sequence. The second one is concerned about the derivation of a
Pohozaev type identity and its application. In subsection 3 we
perform the Blow-up analysis to get either the zero function or a
{\em standard} bubble in the limit. In the subsection 4, we deal
with the capacity estimates to get an upper bound. And in the
subsection 5, we construct test functions . In the last subsection
we show how to reach a contradiction. The last Section is concerned
about the sketch of the proof of Theorem\;$\ref{eq:theorem2}$.
\vspace{10pt}

\noindent {\bf Acknowledgements}

\noindent The second author has been supported by M.U.R.S.T within
the PRIN 2004 Variational methods and nonlinear differential
equations.

\section{Notations and Preliminaries}
In this brief section we collect some useful notations, and state a
lemma giving the existence of the Green functions of \;$\D^2_g$\;and
of the operator \;$P^{4}_{g}$\; with the asymptotics near the
singularity. We also give a version of Adams inequality on the a
manifold due to Fontana and an analogue of the well-known
Moser-Trudinger inequality for the operator \;$P^{4}_{g}$\;when it
is non-negative.

In the following, \;$B_r(x)$\; stands for the metric ball of radius
\;$r$\;and center \;$x$\;in
 \;$M$, $B^r(p)$\; and stands for the Euclidean ball of center \;$p$\; and radius \;$r$. We
 also
 denote with \;$d_{g}(x,y)$\; the metric distance between two points \;$x$\;and \;$y$\;
 of \;$M$.\;$H^{2}(M)$\; stands for the usual  Sobolev space of functions on \;$M$, i.e
 functions
 which are in \;$L^2$\;together with their first and second derivatives. $W^{2,q}(M)$\;denotes
 the usual Sobolev spaces of functions which are in \;$L^q(M)$\;with their first and second
 derivatives. Large positive constants are always denoted by \;$C$,\;and the value of\;$C$\;is
 allowed to vary from formula to formula and also within the same line.\;$M^2$\;stands for the
 cartesian product \;$M\times M$, while \;$Diag(M)$\; is the diagonal of \;$M^2$. Given a
 function \;$u\in L^1(M)$,\;$\bar u$\; denotes its
average on \;$M$,\;that is \;$\bar u=\left(
Vol_{g}(M)\right)^{-1}\int_{M} u(x)dV_{g}(x)$
where \;$Vol_{g}(M)=\int_{M}dV_{g}$.\\
$A_{k}=o_{k}(1)$\; means that \;$A_{k}\rightarrow 0$\;\;as the
integer \;\;$k\rightarrow +\infty$.\\
$A_{\delta}=o_{\delta}(1)$\; means that \;$A_{\delta}\rightarrow
0$\;\;as the real
number \;\;$\delta\longrightarrow 0$.\\
$A_{k,\delta}=o_{k,\delta}(1)$\; means that
\;$A_{k,\delta}\rightarrow 0$\;\;
as \;$k\rightarrow+\infty$\;first and after the real number \;\;$\delta\longrightarrow 0$.\\
$A_{k}=O(B_{k})$\; means that\;$A_{k}\leq C B_{k}$\;\;for some fixed constant \;$C$. \\
$inj_{g}(M)$\;stands for the injectivity radius of \;$M$.\\
$dV_g$\;denotes the Riemannian measure associated to the metric\;$g$.\\
$dS_g$\;stands for the surface measure associated to \;$g$.\\
Given a metric \;$g$\; on \;$M$, and\;$x\in M$,\;$|g(x)|$, stands
for determinant of the matrix with entries \;$g_{i,j}(x)$\;where
\;$g_{i,j}(x)$\;are the components of \;$g(x)$\;in
some system of coordinates.\\
$\D_0$\;stands for the Euclidean Laplacian and \;$\D_g$\;the
Laplace-Beltrami with respect
to the background metric \;$g$.\\

As mentioned before we begin by stating a lemma giving  the
existence of the Green function of
\;$\D^2_g$\;and\;$P^{4}_{g}$,\;and their asymptotics near the
singularities.
\begin{lem}\label{eq:greenlap}
We have that the Green function \;$F(x,y)$\;of \;$\D^2_g$\; exists in the following sense :\\
a) For all functions \;$u\in C^{2}(M)$, we have
\begin{equation*}
u(x)-\bar u=\int_{M}F(x,y)\D^2_{g}u(y)dV_{g}(y)\;\;\;\;\;\;x\neq
y\in M
\end{equation*}
b)
\begin{equation*}
F(x,y)=H(x,y)+K(x,y)
\end{equation*}
is smooth on \;$M^2\setminus Diag(M^2)$,\;$K$\;extends to a
\;$C^{1+\alpha}$\;function on \;$M^2$ \;and
\begin{equation*}
H(x,y)=\frac{1}{8\pi^2}f(r)\log\frac{1}{r}
\end{equation*}
where, \;$r=d_{g}(x,y)$\;is the geodesic distance from \;$x$\;to
\;$y$;\;$f(r)$\;is a \;$C^{\infty}$\;positive  decreasing
function,\;$f(r)=1$\;in a neighborhood of
\;$r=0$\;and\;$f(r)=0$\;for\;$r\geq inj_{g}(M)$. Moreover we have
that the following estimates holds
\begin{equation*}
|\nabla_g F(x,y)|\leq C\frac{1}{d_g(x,y)}\;\;\;|\nabla^2_g F(x,y)
|\leq C\frac{1}{d_g(x,y)^2}.
\end{equation*}
\end{lem}
\begin{proof}
For the proof see \cite{cy} and the proof of Lemma 2.3 in
\cite{mal}.
\end{proof}
\begin{lem}\label{eq:greenpan}
Suppose \;$Ker P^{4}_{g}\simeq\R$. Then the Green function
\;$Q(x,y)$\;of \;$P^{4}_{g}$\;
exists in the following sense :\\
a) For all functions \;$u\in C^{2}(M)$, we have
\begin{equation*}
u(x)-\bar u=\int_{M}Q(x,y)P^{4}_{g}u(y)dV_{g}(y)\;\;\;\;\;\;x\neq
y\in M
\end{equation*}
b)
\begin{equation*}
Q(x,y)=H_0(x,y)+K_0(x,y)
\end{equation*}
is smooth on \;$M^2\setminus Diag(M^2)$,\;$K$\;extends to a
\;$C^{2+\alpha}$\;function on \;$M^2$ \;and
\begin{equation*}
H(x,y)=\frac{1}{8\pi^2}f(r)\log\frac{1}{r}
\end{equation*}
where, \;$r=d_{g}(x,y)$\;is the geodesic distance from \;$x$\;to
\;$y$;\;$f(r)$\;is a \;$C^{\infty}$\;positive  decreasing
function,\;$f(r)=1$\;in a neighborhood of
\;$r=0$\;and\;$f(r)=0$\;for\;$r\geq inj_{g}(M)$.
\end{lem}
\begin{proof}
For the proof see Lemma 2.1 in \cite{nd}.
\end{proof}

Next we state a Theorem due to Fontana\cite{fl}.
\begin{teo}\label{eq:font}(\cite{fl})
There exists a constant \;$C=C(M)>0$\;such that the following holds
\begin{equation*}
\int_{M}e^{32\pi^2u^2}dV_g\leq C\;\;\;\text{for all}\;\;u\in
H^2(M)\;\;\text{such that} \int_{M}|\D_g^2u|dV_g=1.
\end{equation*}
Moreover this constant is optimal in the sense that if we replace it
by any \;$\alpha$\; bigger  then the integral can be maken as large
as we want.
\end{teo}
Next we state a Moser-Trudinger type inequality corresponding to
\;$P^4_g$\;when it is non-negative. The proof can be found in
\cite{nd} where it is proven for every \;$P^n_g$\; (where
\;$P^n_g$\;stands for higher order Paneitz operator).
\begin{prop}
Suppose that \;$P^4_g$\;is non-negative and that \;$ker P^4_g=\R$,
then there exists a constant \;$C=C(M)>0$\;such that
\begin{equation*}
\int_{M}e^{32\pi^2u^2}dV_g\leq C\;\;\;\text{for all}\;\;u\in
H^2(M)\;\;\text{such that} \;\;\left<P^4_gu,u\right>=1.
\end{equation*}
\end{prop}
%%%%%%%%%%%%%%%%%%%%%%%%%%%%%%%%%%%%%%%%%%%%%%%%%%%%%%%%%%%%%%%%%%%%%%%%%%%%%%%%%%%%%

\section{Proof of Theorem\;$\ref{eq:theorem1}$}
\begin{lem}\label{eq:const}
Let \;$\alpha_k$\;be an increasing sequence converging to
\;$32\pi^2$. Then for every \;$k$\; there exists
\;$u_k\in\mathcal{H}_{1}$ \;such that
\begin{equation*}
\int_{M}e^{\alpha_ku_k^2}dV_g=\sup_{u\in\mathcal{H}_{1}}\int_{M}e^{\alpha_ku^2}dV_g.
\end{equation*}
Moreover \;$u_k$\;satisfies the following equation
\begin{equation}\label{eq:edp}
\D_g^2u_k=\frac{1}{\lambda_{k}}u_ke^{\alpha_ku_k^2}-\gamma_k
\end{equation}
where
\begin{equation*}
\lambda_{k}=\int_{M}u_k^2e^{\alpha_ku_k^2}dV_g
\end{equation*}
and
\begin{equation*}
\gamma_k=\frac{1}{\lambda_kVol_g(M)}\int_{M}u_ke^{\alpha_ku_k^2}dV_g.
\end{equation*}
Moreover we have \;$u_k\in C^{\infty}(M)$.
\end{lem}
\begin{proof}
First of all using the inequality in Theorem\;$\ref{eq:font}$, one
can check easily that the functional
\begin{equation*}
I_k(u)=\int_{M}e^{\alpha_ku^2}dV_g;
\end{equation*}
is weakly continous. Hence using Direct Methods of the Calculus of
Variations  we get the existence of maximizer say \;$u_k$. On the
other hand using the Lagrange multiplier rule one get the equation
(\ref{eq:edp}). Moreover integrating the equation  (\ref{eq:edp})
and after multipling it by \;$u_k$\;and integrating again, we get
the value
 of \;$\gamma_k$\;and \;$\lambda_k$\; respectively. Moreover using standard elliptic
regularity we get that \;$u_k\in C^{\infty}(M)$. Hence the Lemma is
proved.
\end{proof}
Now we are ready to give the proof of Theorem\;$\ref{eq:theorem1}$.
>From now on we suppose by contradiction that
Theorem\;$\ref{eq:theorem1}$\;does not hold. Hence from the same
considerations as in the Introduction we
have that :\\
1)\\
\begin{equation}\label{eq:blowup}
\forall
\alpha>32\pi^2\;\;\;\lim_{k\rightarrow+\infty}\int_{M}e^{\alpha
u_k^2}dV_g\rightarrow+\infty
\end{equation}
2)\\
\begin{equation*}
c_k=\max_{M}|u_k|=|u_k|(x_k)\rightarrow +\infty
\end{equation*}
We will divide the reminder of the proof into six subsections.
\subsection{Concentration behavior and profile of\;$u_k$}
This subsection is concerned about two main ingredients. The first
one is the study of the concentration  phenomenon of the energy
corresponding to\;$u_k$. The second one is the description of the
profile of \;$\beta_ku_k$\; as \;$k\rightarrow+\infty$, where
\;$\beta_k$\;is given by the relation
\begin{equation*}
1/\beta_k=\int_M\frac{|u_k|}{\lambda_k}e^{\alpha_ku_k^2}dV_g.
\end{equation*}
We start by giving an energy concentration lemma which is inspired
from P.L.Lions'work.
\begin{lem}\label{eq:conc-comp}
\;$u_k$\;verifies :
\begin{equation*}
u_k\rightharpoonup 0\;\;\text{in}\;\;H^2(M);
\end{equation*}
and
\begin{equation*}
|\D_g u_k|^2 \rightharpoonup \delta_{x_0}
\end{equation*}
for some \;$x_0\in M$.
\end{lem}
\begin{proof}
First of all from the fact that \;$u_k\in \mathcal{H}_{1}$\; we can
assume  without loss of generality that
\begin{equation}\label{uk}
u_k\rightharpoonup u_{0}\;\;\text{in}\;\;H^2(M).
\end{equation}
Now let us show that \;$u_0=0$.\\
We have the trivial identity
\begin{equation*}
\int_{M}|\D_g(u_k-u_0)|^2dV_g=\int_{M}|\D_gu_k|^2dV_g
+\int_{M}|\D_gu_0|^2dV_g-2\int_{M} \D_gu_k\D_gu_0dV_g.
\end{equation*}
Hence using the fact that \;$\int_{M}|\D_gu_k|^2dV_g=1$\; we derive
\begin{equation*}
\int_{M}|\D_g(u_k-u_0)|^2dV_g=1
+\int_{M}|\D_gu_0|^2dV_g-2\int_{M}\D_gu_k\D_gu_0dV_g
\end{equation*}
So using  \;(\ref{uk})\; we get
\begin{equation*}
\lim_{k\rightarrow
0}\int_{M}|\D_g(u_k-u_0)|^2dV_g1-\int_{M}\D_gu_{0}\D_gu_0dV_g
\end{equation*}
Now suppose that \;$u_0\neq 0$\; and let us argue for a
contradiction. Then there exists some \;$\beta<1$\;such that for
\;$k$\;large enough the following holds
\begin{equation*}
\int_{M}|\D_g(u_k-u_0)|^2dV_g<\beta.
\end{equation*}
Hence using Fontana's result see Theorem\;$\ref{eq:font}$\;we obtain
that
\begin{equation*}
\int_{M}e^{\alpha_1 (u_k-u_{0})^2}dV_g\leq C \hbox{ for some }
\alpha_1>32\pi^2.
\end{equation*}
Now using Cauchy inequality one can check easily that
\begin{equation*}
\int_{M}e^{\alpha_2 u_k^2}dV_g\leq C \hbox{ for some }
\alpha_2>32\pi^2.
\end{equation*}
Hence reaching a contradiction to\;$\eqref{eq:blowup}$.\\
On the other hand without lost of generality we can assume that
\begin{equation*}
|\D_g u_k|dV_g\rightharpoonup \mu.
\end{equation*}
Now suppose \;$\mu\neq \delta_{p}$\;for every \;$p\in M$\;and let us
argue for a contradiction to \;$\eqref{eq:blowup}$\; again. First of
all let us take a cut-off function \;$\eta\in
C^{\infty}_{0}(B_{\delta}(x))$\;,\;$\eta=1$\;on
\;$B_{x}(\frac{\delta}{2})$\;
where \;$x$\;is a fixed point in \;$M$\; and \;$\delta$\; a fixed positive and small number.\\
We have that
\begin{equation*}
\limsup_{k\rightarrow+\infty}\int_{B_{\delta}(x)}|\D_gu_k|^2dV_g<1.
\end{equation*}
Now working in a normal coordinate system around \;$x$\;and using
standard elliptic regularity theory we get
\begin{equation*}
\int_{B^{\delta}(\tilde x)}|\D_0\widetilde{ \eta u_k}|^2dV_g\leq
(1+o_{\delta}(1)) \int_{B_{\delta}(x)}|\D_gu_k|^2dV_g;
\end{equation*}
where \;$\tilde x$\;is the  point corresponding to \;$x$\; in
\;$\R^4$\; and \;$\widetilde{ \eta u_k}$\; the expression of\;$ \eta
u_k$\;on the normal coordinate system. Hence for \;$\delta \;$small
we get
\begin{equation*}
\int_{B^{\delta}(\tilde x)}|\D_0\widetilde{ \eta u_k}|^2dV_g<1
\end{equation*}
Thus using the Adams result see \cite{adr} we have that
\begin{equation*}
\int_{B^{\delta}(\tilde x)}e^{\tilde \alpha(\widetilde{ \eta
u_k)^2}}dx\leq C \hbox{ for some } \tilde{\alpha}>32\pi^2.
\end{equation*}
Hence using a covering argument we infer that
\begin{equation*}
\int_{M}e^{\bar \alpha u_k^2}dV_g\leq C  \hbox{for some} \bar
\alpha>32\pi^2,
\end{equation*}
so reaching a contradiction. Hence the Lemma is proved.
\end{proof}
\begin{lem}\label{eq:asymp}
We have the following hold:
\begin{equation*}
\lim_{k\rightarrow
+\infty}\lambda_{k}=+\infty,\;\;\;\;\lim_{k\rightarrow
+\infty}\gamma_k=0.
\end{equation*}
 \end{lem}
\begin{proof}
Let \;$N>0$\;be large enougth. By using the definition of
\;$\lambda_{k}$\;we have that
\begin{equation*}
\lambda_{k}=\int_{M}u_k^2e^{\alpha_ku_k^2}dV_g\geq
N^2\int_{\{u_k\geq N\}}e^{\alpha_ku_k^2}dV_g
=N^2(\int_{M}e^{\alpha_ku_k^2}dV_g-\int_{\{u_k\leq N\}}
e^{\alpha_ku_k^2}dV_g).
\end{equation*}
On the other hand
\begin{equation*}
\lim_{k\rightarrow
+\infty}\left(\int_{M}e^{\alpha_ku_k^2}dV_g-\int_{\{u_k\leq N\}}
e^{\alpha_ku_k^2}dV_g\right)=\lim_{k\rightarrow+\infty}\int_{M}e^{\alpha_ku_k^2}dV_g-Vol_g(M).
\end{equation*}
Hence using the fact that
\begin{equation*}
\lim_{k\rightarrow+\infty}\int_{M}e^{\alpha_ku_k^2}dV_g=\sup_{u\in\mathcal{H}_{1}}
\int_{M}e^{32\pi^2u^2} dV_g>Vol_g(M)
\end{equation*}
we have that 1) holds. Now we  prove 2). using the definition of
\;$\gamma_k$\;, we get
\begin{equation*}
|\gamma_k|\leq
\frac{N}{\lambda_{k}}Ne^{32\pi^2N^2}+\frac{1}{Vol_g(M)}\frac{1}{N}.
\end{equation*}
Hence by using point 1 and letting \;$k\rightarrow +\infty$\;and
after \;$N\rightarrow +\infty$\;we get point 2. So the Lemma is
proved.
\end{proof}
Next let us set
\begin{equation*}
\tau_k=\int_M\frac{\beta_ku_k}{\lambda_k}e^{\alpha_ku_k^2}.
\end{equation*}
One can check easily the following
\begin{lem}\label{eq:asymli}
With the definition above we have that $0\leq \beta_{k}\leq
c_{k}$,\;\;$|\tau_{k}|\leq 1$\;and \;$\beta_{k}\gamma_k$\;is
bounded. Moreover up to a subsequence and up to changing
\;$u_{k}$\;to \;$-u_{k}$
\begin{equation*}
\tau_{k}\rightarrow \tau \geq 0.
\end{equation*}
\end{lem}

The next Lemma gives some Lebesgue estimates on Ball in terms of the
radius with constant independent of the ball. As a corollary we get
the profile of \;$\beta_ku_k$\;as \; $k\rightarrow +\infty$.
\begin{lem}\label{eq:lpest}
There are constants $C_1(p)$,and\; $C_2(p)$\; depending only on
\;$p$\; and \;$M$\; such that, for $r$ sufficiently small and for
any $x\in M$ there holds
\begin{equation*}
\int_{B_r(x)}|\nabla_g^2\beta_{k}u_{k}|^pdV_g\leq C_2(p)r^{4-2p};
\end{equation*}
and
\begin{equation*}
\int_{B_r(x)}|\nabla_g\beta_{k} u_{k}|^pdV_g\leq C_1(p)r^{4-p}
\end{equation*}
where, respectively, $p<2$, and $p<4$.
\end{lem}
\begin{proof}
First of all using the Green representation formula we have
\begin{equation*}
u_{k}(x)=\int_{M}F(x,y)\D_g^2u_{k}dV_g(y)\;\;\;\forall x\in M.
\end{equation*}
Hence using the equation we get
\begin{equation*}
u_{k}(x)=\int_{M}F(x,y)\left(\frac{1}{\lambda_{k}}u_{k}e^{\alpha_ku_{k}^2}\right)dV_g(y)-
\int_{M}F(x,y)\gamma_kdV_g(y).
\end{equation*}
Now by differentiating with respect to \;$x$\;for every
\;$m=1,2$\;we have that
\begin{equation*}
|\nabla_g^m u_{k}(x)|\leq
\int_{M}|\nabla_g^{m}F(x,y)|\left(\frac{1}{\lambda_{k}}\right)
|u_{k}|e^{\alpha_ku_{k}^2}dV_g(y)+\int_{M}|\nabla_g^m
F(x,y)|\left|\gamma_k\right|.
\end{equation*}
Hence we get
\begin{equation*}
|\nabla_g^m (\beta_{k}u_{k}(x))|\leq
\int_{M}|\nabla_g^{m}F(x,y)|\beta_{k}\left(\frac{1}
{\lambda_{k}}\right)|u_{k}|e^{\alpha_ku_{k}^2}dV_g(y)+\int_{M}|\nabla_g^m
F(x,y)|\beta_{k} \left|\gamma_k\right|.
\end{equation*}
Taking the \;$p$-th power in both side of the inequality and using
the basic inequality
\begin{equation*}
(a+b)^p\leq 2^{p-1}(a^p+b^p) \;\;\;\text{for}\;\;\; a\geq
0\;\;\;\text{and}\;\;\ b\geq 0
\end{equation*}
we obtain
\begin{equation*}
\begin{split}
|\nabla_g^m (\beta_{k}u_{k}(x))|^p\leq 2^{p-1}\left[
\int_{M}|\nabla_g^{m}F(x,y)|\beta_{m}
\left(\frac{1}{\lambda_{k}}\right)|u_{k}|e^{\alpha_ku_{k}^2}dV_g(y)\right]^p\\+2^{p-1}
\left[\int_{M}|\nabla_g^m
F(x,y)|\beta_{k}\left|\gamma_k\right|\right]^p
\end{split}
\end{equation*}
Now integrating both sides of the inequality we obtain
\begin{equation*}
\begin{split}
\int_{B_{r}(x)} |\nabla_g^m (\beta_{k}u_{k}(z))|dV_g(z)\leq
2^{p-1}\int_{B_r(x)}\left[
\int_{M}|\nabla_g^{m}F(z,y)|\beta_{k}\left(\frac{1}{\lambda_{k}}\right)
|u_{k}|e^{\alpha_ku_{k}^2}dV_g(y)\right]^pdV_g(z)\\+2^{p-1}\int_{B_{r}(x)}
\left[\int_{M}|\nabla_g^m
F(z,y)|\beta_{k}\left|\gamma_k\right|\right]^pdV_g(z).
\end{split}
\end{equation*}
First let us estimate the second term in the right hand side of the
inequality
\begin{equation*}
\int_{B_r(x)}\left[\int_{M}|\nabla_g^m
F(z,y)|\beta_{k}\left|\gamma_k\right|\right]^pdV_g(z) \leq C
\int_{B_r(x)}\sup_{y\in M}\frac{1}{d_{g}(z,y)^{pm}}dV_g(z)\leq
C(M)r^{4-mp}
\end{equation*}
Thanks to the fact that \;$\beta_k\gamma_k$\;is bounded, to the
asymptotics of the Green function and to Jensen's inequality. Now
let us estimates the second term. First of all we define the
following auxiliary measure
\begin{equation*}
m_{k}=\beta_{k}\left(\frac{1}{\lambda_{k}}\right)|u_{k}|e^{\alpha_ku_{k}^2}dV_g
\end{equation*}
We have that \;$m_{k}$\;is a probability measure. On the other hand
we can write
\begin{equation}\label{jen1}
\begin{split}
\int_{B_r(x)}\left[
\int_{M}|\nabla_g^{m}F(z,y)|\beta_{k}\left(\frac{1}{\lambda_{k}}
\right)|u_{k}|e^{\alpha_ku_{k}^2}dV_g(y)\right]^pdV_g(z)\\=\int_{B_r(x)}\left[
\int_{M} |\nabla_g^{m}F(z,y)|dm_{k}(y)\right]^pdV_g(z).
\end{split}
\end{equation}
Now by using Jensen's inequality we have that
\begin{equation*}
\left[ \int_{M}|\nabla_g^{m}F(z,y)|dm_{k}(y)\right]^p\leq \left[
\int_{M}|\nabla_g^{m} F(z,y)|^pdm_{k}(y)\right]
\end{equation*}
Thus with the (\ref{jen1}) we have that
\begin{equation*}
\begin{split}
\int_{B_r(x)}\left[
\int_{M}|\nabla_g^{m}F(z,y)|\beta_{k}\left(\frac{1}{\lambda_{k}}
\right)|u_{k}|e^{\alpha_ku_{k}^2}dV_g(y)\right]^pdV_g(z)
\leq\\\int_{B_r(x)}\left[
\int_{M}|\nabla_g^{m}F(z,y)|^pdm_{k}(y)\right]dV_g(z).
\end{split}
\end{equation*}
Now by using again the same argument as in the first term we obtain
\begin{equation*}
\int_{B_r(x)}\left[
\int_{M}|\nabla_g^{m}F(z,y)|^pdm_{k}(y)\right]dV_g(z)\leq
C(M)r^{4-mp}.
\end{equation*}
Hence the Lemma is proved.
\end{proof}
Next we give a corollary of this Lemma.
\begin{cor}\label{eq:profile}
We have \;$\beta_{k}u_{k}\rightharpoonup G$\;\;$W^{2,p}(M)$ \;for
\;$p\in (1,2)$, \;$\beta_{k}u_{k} \rightharpoonup G$ smoothly in
\;$M\setminus B_{\delta}(x_0)$\; where \;$\delta $\;is small and
\;$G$\;satisfies
\begin{equation*}
\left\{
\begin{array}{ll}
 \D_g^2 G=\tau(\delta_{x_0}-\frac{1}{Vol_g(M)})\;\;\text{in}\;\;M;&\\
 \overline G=0
\end{array}
\right.
\end{equation*}
Moreover
\begin{equation*}
G(x)=\frac{\tau}{8\pi^2}\log\frac{1}{r}+\tau S(x)
\end{equation*}
with \;$r=d_g(x,x_0)$.\;$S=S_0+S_1(x)$\;,\;$S_0=S(x_0)$\;and \;$S\in
W^{2,q}(M)$\;for every \;$q\geq 1$.
\end{cor}
\begin{proof}
By Lemma\;$\ref{eq:lpest}$\; we have that
\begin{equation*}
\beta_{k}u_{k}\rightharpoonup G\;\;W^{2,p}(M)\;\;p\in (1,2)
\end{equation*}
On the other hand using Lemma\;$\ref{eq:conc-comp}$\; we
get\;$e^{\alpha_ku_{k}^2}$\;is bounded in \;$L^p(M\setminus
B_{\delta}(x_0))$. Hence the standard elliptic regularity implies
that
\begin{equation}\label{bu}
\beta_{k}u_{k}\rightarrow G\;\;\text{smoothly in }\;M\setminus
B_{\delta}(x_0).
\end{equation}
So to end the proof of the proposition we need only to show that
\begin{equation}\label{claim}
\frac{\beta_{k}}{\lambda_{k}}u_{k}e^{\alpha_ku_{k}^2}\rightharpoonup
\tau \delta_{x_0}.
\end{equation}
To do this let us take \;$\varphi\in C^{\infty}(M)$\; then we have
\begin{equation*}
\int_{M}\varphi
\frac{\beta_{k}}{\lambda_{k}}u_{k}e^{\alpha_ku_{k}^2}dV_g=\int_{M\setminus
B_{\delta}(x_0)}\varphi
\frac{\beta_{k}}{\lambda_{k}}u_{k}e^{\alpha_ku_{k}^2}dV_g+\int_{
B_{\delta}(x_0)}\varphi
\frac{\beta_{k}}{\lambda_{k}}u_{k}e^{\alpha_ku_{k}^2}dV_g
\end{equation*}
Using (\ref{bu})  we have that
\begin{equation*}
\int_{M\setminus B_{\delta}(x_0)}\varphi
\frac{\beta_{k}}{\lambda_{k}}u_{k}e^{\alpha_ku_{k}^2}
dV_g=O(\frac{1}{\lambda_{k}}).
\end{equation*}
On the other hand, we can write inside the ball \;$B_{\delta}(x_0)$
\begin{equation*}
\begin{array}{lll}
\int_{ B_{\delta}(x_0)}\varphi
\frac{\beta_{k}}{\lambda_{k}}u_{k}e^{\alpha_ku_{k}^2}dV_g
&=&(\varphi(x_0)+o_{\delta}(1))\dint_{ B_{\delta}(x_0)}
\frac{\beta_{k}}{\lambda_{k}}u_{k}
e^{\alpha_ku_{k}^2}dV_g\\[\mv]
&=&(\varphi(x_0)+o_{\delta}(1))\left(\tau-\dint_{M\setminus
B_{\delta}(x_0)} \frac{\beta_{k}}
 {\lambda_{k}}u_{k}e^{\alpha_ku_{k}^2}dV_g\right)
\end{array}
\end{equation*}
Now using again \;$(\ref{bu})$\; we derive
\begin{equation*}
\int_{M\setminus B_{\delta}(x_0)}
\frac{\beta_{k}}{\lambda_{k}}u_{k}e^{\alpha_ku_{k}^2}
=O(\frac{1}{\lambda_k}).
\end{equation*}
Hence  we arrive to
\begin{equation*}
\int_{ B_{\delta}(x_0)}\varphi
\frac{\beta_{k}}{\lambda_{k}}u_{k}e^{\alpha_ku_{k}^2}dV_g =\tau
\varphi(x_0)+o_{k,\delta}(1).
\end{equation*}
Thus we get
\begin{equation*}
\int_{M}\varphi
\frac{\beta_{k}}{\lambda_{k}}u_{k}e^{\alpha_ku_{k}^2}dV_g=O(\frac{1}
{\lambda_{k}})+\tau \varphi(x_0)+o_{k,\delta}(1).
\end{equation*}
Hence from Lemma\;$\ref{eq:asymp}$\; we conclude the proof of claim
(\ref{claim}) )and of the Corollary too.
\end{proof}
\subsection{Pohozaev type identity and application}
As it is already said in the introduction this subsection deals with
the derivation of a Pohozaev type identity. And as corollary we give
the limit of \;$\int_{M}e^{\alpha_ku_k^2}dV_g$\;in terms of
\;$Vol_g(M)$,\;$\lambda_k$,\;$\beta_k$\; and \;$\tau$
\begin{lem}\label{eq:P}
Setting \; $U_{k}=\D_g u_{k}$\;we have the following holds
\begin{equation*}
\begin{split}
    -\frac{2}{\alpha_k\lambda_k}
      \int_{B_\delta(x_k)}e^{\alpha_ku_k^2}dV_g=-\frac{\delta}{2}
       \int_{\partial B_\delta(x_k)}
       U_k^2dS_g-\delta\int_{\partial B_\delta(x_k)}\nabla_gu_k
       \nabla_gU_kdV_g+
       2\int_{\partial B_\delta(x_k)}U_k\frac{\partial u_k}
       {\partial r}\\
    +2\delta\int_{\partial {B_\delta(x_k)}}\frac{\partial U_k}{\partial r}
          \frac{\partial u_k}{\partial r}dS_g
          +\int_{B_\delta(x_k)}O(r^2)\nabla_g u_k \nabla_g U_k dV_g\\
    +\int_{B_\delta(x_k)}O(r^2)U_k^2dV_g+
      \int_{B_\delta(x_k)}
         e^{\alpha_ku_k^2}O(r^2)dV_g\
     -\frac{\delta}{2\lambda_k\alpha_k}\int_{\partial B_{\delta}(x_k)}
     e^{\alpha_ku_k^2}dV_g+O(\frac{\delta}{\beta_k^2}).
      \end{split}
\end{equation*}
where\;$\delta$\;is small and fixed real number.
\end{lem}
\begin{proof}
The proof relies on the divergence formula and the asymptotics of
the metric g in normal
coordinates around \;$x_k$.\\
By the definition of \;$U_{k}$ we have that
\begin{equation*}
\left\{\begin{array}{l}
            \Delta_g u_k=U_k\\[1.5ex]
            \Delta_g U_k=\frac{u_k}{\lambda_k}e^{\alpha_ku_k^2}-\gamma_k.
          \end{array}\right.
      \end{equation*}

 The first issue is to compute \; $\int_{B_\delta(x_k)}r\frac{\partial U_k}{\partial r}
\Delta_g u_k$\; in two different ways, where $r(x)=d_g(x,x_k)$.\\
On one side we obtain
$$\begin{array}{ll}
      \dint_{B_\delta(x_k)}r\frac{\partial U_k}{\partial r}\Delta_g u_kdV_g
        &=-\dint_{B_\delta(x_k)}(\nabla_g U_k\nabla_g u_k+r\frac{\partial\nabla_g U_k}
           {\partial r}\nabla_g u_k)dV_g
          +\dint_{\partial {B_\delta(x_k)}}r\frac{\partial U_k}{\partial r}
          \frac{\partial u_k}{\partial r}dS_g.
   \end{array}$$
On the other side we get
$$\begin{array}{ll}
      \dint_{B_\delta(x_k)}r\frac{\partial U_k}{\partial r}\Delta_g u_kdV_g
        &=\dint_{{B_\delta(x_k)}}r\frac{\partial U_k}{\partial r}U_kdV_g\\[\mv]
        &=\dint_0^\delta2\pi^2 \dint_{\partial B_r(x_k)}\frac{\partial U_k}{\partial r}
          U_k\sqrt{|g|}r^4dSdr\\[\mv]
        &=\frac{\delta}{2}\dint_{\partial B_\delta(x_k)}U_k^2dS_g-2\dint_{B_\delta(x_k)}
          U_k^2(1+O(r^2))dV_g.
  \end{array}$$

Thus we have
$$\begin{array}{lll}
  \frac{\delta}{2}\dint_{\partial B_\delta(x_k)}U_k^2dS_g-2\dint_{B_\delta(x_k)}U_k^2dV_g&=&
       -\dint_{B_\delta(x_k)}(\nabla_g U_k\nabla_g u_k+r\frac{\partial\nabla_g U_k}
       {\partial r}\nabla_g u_k)dV_g\\[\mv]
  &&+\dint_{\partial {B_\delta(x_k)}}r\frac{\partial U_k}{\partial r}
          \frac{\partial v_k}{\partial r}dS_g
    +\dint_{B_\delta(x_k)}O(r^2)U_k^2dV_g
  \end{array}$$
In the same way we obtain
$$\begin{array}{l}
    \frac{\delta}{2\lambda_k\alpha_k}\dint_{\partial B_\delta(x_k)}
        e^{\alpha_ku_k^2}dS_g-\frac{2}{\lambda_k\alpha_k}\dint_{B_\delta(x_k)}
         e^{\alpha_ku_k^2}(1+O(r^2))dV_g\\[\mv]
   \s\s\s=-\dint_{B_\delta(x_k)}(\nabla_g U_k\nabla_g u_k+r\frac{\partial\nabla_g u_k}
   {\partial r}
     \nabla_g U_k)dV_g+\dint_{\partial {B_\delta(x_k)}}r\frac{\partial U_k}{\partial r}
          \frac{\partial u_k}{\partial r}dS_g +O(\frac{\delta}{\beta_k^2}).
  \end{array}$$
Hence by summing this two last lines we arrive to
\begin{equation}\label{eq:1}
\begin{array}{l}
     \frac{\delta}{2\lambda_k\alpha_k}\dint_{\partial B_\delta(x_k)}
        e^{\alpha_ku_k^2}dS_g-\frac{2}{\lambda_k\alpha_k}\dint_{B_\delta(x_k)}
         e^{\alpha_ku_k^2}dV_g+
   \frac{\delta}{2}\dint_{\partial B_\delta(x_k)}U_k^2dS_g-2\dint_{B_\delta(x_k)}U_k^2dV_g\\[\mv]
   \s\s=-\dint_{B_\delta(x_k)}(2\nabla_g U_k\nabla_g u_k+r\frac{\partial}{\partial r}
    \nabla_g u_k \nabla_g U_k)dV_g+2\dint_{\partial {B_\delta(x_k)}}r\frac{\partial U_k}
    {\partial r}
          \frac{\partial u_k}{\partial r}dS_g\\[\mv]
   \s\s\s\s+\dint_{B_\delta(x_k)}O(r^2)U_k^2dV_g+
      \dint_{B_\delta(x_k)}
         e^{\alpha_ku_k^2}O(r^2)dV_g+O(\frac{\delta}{\beta_k^2}).
  \end{array}
  \end{equation}
On the other hand  using  the same method one can check easily that
\begin{equation}\label{eq:2}
\begin{array}{lll}
\dint_{B_\delta(x_k)}r\frac{\partial}{\partial r}
    \nabla_g u_k \nabla_g U_kdV_g&=&\delta\dint_{\partial B_\delta(x_k)}
     \nabla_g u_k\nabla_gU_k dV_g-4\dint_{B_\delta(x_k)}\nabla_gu_k\nabla_gU_kdV_g\\[\mv]
     &&+\dint_{B_\delta(x_k)}O(r^2)\nabla_g u_k \nabla_g U_k dV_g
  \end{array}
  \end{equation}
and
\begin{equation}\label{eq:3}
\begin{array}{lll}
   \dint_{B_\delta(x_k)}\nabla_g U_k\nabla_g u_kdV_g&=&
       -\dint_{B_\delta(x_k)}U_k\Delta_g u_kdV_g+
        \dint_{\partial B_\delta(x_k)}U_k\frac{\partial u_k}{\partial r}dS_g\\[\mv]
    &=&-\dint_{B_\delta(x_k)}U_k^2dV_g+
       \dint_{\partial B_\delta(x_k)}U_k\frac{\partial u_k}
       {\partial r}dS_g,
  \end{array}
  \end{equation}
So using\;$\eqref{eq:1}$,$\eqref{eq:2}$\;and \;$\eqref{eq:3}$\;we
arrive to
$$\begin{array}{lll}
    -\frac{2}{\alpha_k\lambda_k}
      \dint_{B_\delta(x_k)}e^{\alpha_ku_k^2}dV_g&=&-\frac{\delta}{2}
       \dint_{\partial B_\delta(x_k)}
       U_k^2dS_g-\delta\dint_{\partial B_\delta(x_k)}\nabla_gu_k
       \nabla_gU_kdV_g+
       2\dint_{\partial B_\delta(x_k)}U_k\frac{\partial u_k}
       {\partial r}\\[\mv]
    &&+2\delta\dint_{\partial {B_\delta(x_k)}}\frac{\partial U_k}{\partial r}
          \frac{\partial u_k}{\partial r}dS_g
          +\dint_{B_\delta(x_k)}O(r^2)\nabla_g u_k \nabla_g U_k dV_g\\[\mv]
    &&+\dint_{B_\delta(x_k)}O(r^2)U_k^2dV_g+
      \dint_{B_\delta(x_k)}
         e^{\alpha_ku_k^2}O(r^2)dV_g\\[\mv]
     &&-\frac{\delta}{2\lambda_k\alpha_k}\dint_{\partial B_{\delta}(x_k)}
     e^{\alpha_ku_k^2}dV_g+O(\frac{\delta}{\beta_k^2}).
      \end{array}$$
Thus the Lemma is proved
\end{proof}
\begin{cor}\label{eq:sup}
 We have that
 \begin{equation*}
 \lim_{k\rightarrow +\infty}
\int_Me^{\alpha_ku_k^2}=Vol_g(M)+\tau^2\lim_{k\rightarrow+\infty}\frac{\lambda_k}{\beta_k^2}.
 \end{equation*}
 Moreover we have that
 \begin{equation*}
 \tau\in (0,1].
 \end{equation*}
\end{cor}
\begin{proof}
First of all we have that the sequence
\;$(\frac{\lambda_k}{\beta_k^2})_k$\;is bounded. Indeed using the
definition of \;$\beta_k$\;we have that
\begin{equation*}
\frac{\lambda_k}{\beta_k^2}=\frac{1}{\lambda_k}(\int_{M}|u_k|e^{\alpha_ku_k^2}dV_g)^2.
\end{equation*}
Hence using Jensen's inequality we obtain
\begin{equation*}
\frac{\lambda_k}{\beta_k^2}\leq\frac{1}{\lambda_k}\int_{M}e^{\alpha_ku_k^2}dV_g\int_{M}
u_k^2e^{\alpha_ku_k^2}dV_g.
\end{equation*}
Thus using the definition of \;$\lambda_k$\;we have that
\begin{equation*}
\frac{\lambda_k}{\beta_k^2}\leq\int_{M}e^{\alpha_ku_k^2}dV_g.
\end{equation*}
On the other hand one can check easily that
\begin{equation*}
\lim_{k\rightarrow+\infty}\int_{M}e^{\alpha_ku_k^2}dV_g=\sup_{u\in\mathcal{H}_{1}}
\int_{M}e^{32\pi^2u^2}dV_g<\infty.
\end{equation*}
Hence we derive that \;$(\frac{\lambda_k}{\beta_k^2})_k$\;is
bounded. So we can suppose
without lost of generality that \;$(\frac{\lambda_k}{\beta_k^2})_k$\;converges.\\
Now from Lemma \ref{eq:P} we have that
\begin{equation*}
\begin{array}{lll}
   \lim\limits_{k\rightarrow +\infty}
      \dint_{B_\delta(x_k)}e^{\alpha_ku_k^2}dV_g&=&16\pi^2\lim\limits_{k\rightarrow+\infty}
      \frac{\lambda_k}{\beta_k^2}(\frac{\delta}{2}
       \dint_{\partial B_\delta(x_k)}
       (\beta_kU_k)^2dS_g \\[\mv]&&+\delta\dint_{\partial B_\delta(x_k)}\nabla_g(\beta_ku_k)
       \nabla_g(\beta_kU_k)dS_g-
       2\dint_{\partial B_\delta(x_k)}(\beta_kU_k)\frac{\partial (\beta_ku_k)}
       {\partial r}\\[\mv]
    &&-2\delta\dint_{\partial {B_\delta(x_k)}}\frac{\partial (\beta_kU_k)}{\partial r}
          \frac{\partial (\beta_ku_k)}{\partial r}dS_g
          +O(\delta)).
      \end{array}
\end{equation*}
So using Lemma \ref{eq:profile} we obtain
\begin{equation*}
\begin{array}{lll}
   \lim\limits_{k\rightarrow +\infty}
      \dint_{B_\delta(x_k)}e^{\alpha_ku_k^2}dV_g&=&16\pi^2\lim\limits_{k\rightarrow+\infty}
      \frac{\lambda_k}{\beta_k^2}(\frac{\delta}{2}
       \dint_{\partial B_\delta(x_0)}
       |\D_g G|^2dS_g \\[\mv]&&+\delta\dint_{\partial B_\delta(x0)}\nabla_gG
       \nabla_g(\D_g G)dS_g-
       2\dint_{\partial B_\delta(x_0)}\D_g G\frac{\partial G}
       {\partial r}\\[\mv]
    &&-2\delta\dint_{\partial {B_\delta(x_0)}}\frac{\partial \D_gG}{\partial r}
          \frac{\partial G}{\partial r}dS_g
          +O(\delta)).
      \end{array}
\end{equation*}
Moreover by trivial calculations we get
\begin{equation*}
\int_{\partial B_\delta(x_0)}
       |\D_g G|^2dS_g =\frac{\tau^2}{8\pi^2\delta}+O(1);
\end{equation*}
\begin{equation*}
\int_{\partial B_\delta(x_0)}\nabla_gG
       \nabla_g(\D_g G)dS_g=-\frac{\tau^2}{8\pi^2\delta}+O(1);
\end{equation*}
\begin{equation*}
\int_{\partial B_\delta(x_0)}\D_g G\frac{\partial G}
       {\partial r}=\frac{\tau^2}{16\pi^2}+O(\delta);
\end{equation*}
and
\begin{equation*}
\int_{\partial {B_\delta(x_0)}}\frac{\partial \D_g G}{\partial r}
          \frac{\partial G}{\partial r}dS_g=-\frac{\tau^2}{8\pi^2\delta}+O(1)
\end{equation*}
Hence with this we obtain
\begin{equation*}
\lim\limits_{k\rightarrow +\infty}
      \int_{B_\delta(x_k)}e^{\alpha_ku_k^2}dV_g=\tau^2\lim\limits_{k\rightarrow+\infty}
      \frac{\lambda_k}{\beta_k^2}+O(\delta).
\end{equation*}
On the other hand  we have that
\begin{equation*}
\int_{M}e^{\alpha_ku_k^2}dV_g=\int_{B_\delta(x_k)}e^{\alpha_ku_k^2}dV_g+
\int_{M\setminus B_\delta(x_k)}e^{\alpha_ku_k^2}dV_g
\end{equation*}
Moreover by Lemma\;$\ref{eq:conc-comp}$\; we have that
\begin{equation*}
\int_{M\setminus
B_\delta(x_k)}e^{\alpha_ku_k^2}dV_g=Vol_g(M)+o_{k,\delta}(1).
\end{equation*}
Thus we derive that
\begin{equation*}
\lim\limits_{k\rightarrow
+\infty}\int_{M}e^{\alpha_ku_k^2}dV_g=Vol_g(M)+
\tau^2\lim\limits_{k\rightarrow+\infty}\frac{\lambda_k}{\beta_k^2}+o_{\delta}(1).
\end{equation*}
Hence letting \;$\delta\rightarrow 0$\;we obtain
\begin{equation*}
\lim\limits_{k\rightarrow
+\infty}\int_{M}e^{\alpha_ku_k^2}dV_g=Vol_g(M)+
\tau^2\lim\limits_{k\rightarrow+\infty}\frac{\lambda_k}{\beta_k^2}.
\end{equation*}
Now suppose \;$\tau=0$\;then we get
\begin{equation*}
\lim\limits_{k\rightarrow
+\infty}\int_{M}e^{\alpha_ku_k^2}dV_g=Vol_g(M).
\end{equation*}
On the other hand we have that
\begin{equation*}
\lim\limits_{k\rightarrow
+\infty}\int_{M}e^{\alpha_ku_k^2}dV_g=\sup_{u\in\mathcal
{H}_1}\int_{M}e^{32\pi^2u^2}dV_g>Vol_g(M);
\end{equation*}
hence a contradiction. Thus \;$\tau\neq 0$\;and the Corollary is
proved.
\end{proof}
\subsection{Blow-up analysis}
In this subsection we perform the Blow-up analysis and show that the
asymptotic profile of \;$u_k$\;is either
the zero function or a standard Bubble.\\
First of all let us introduce some notations.\\
We set
\begin{equation*}
r_k^4=\frac{\lambda_k}{\beta_kc_k}e^{-\alpha_kc_k^2}.
\end{equation*}
Now  for \;$x\in B^{r_k^{-1}\delta}(0)$\;with \;$\delta>0$\;small we
set
\begin{equation*}
w_{k}(x)=2\alpha_k\beta_{k}\left(u_{k}(exp_{x_{k}}(r_{k}x))-c_{k}\right);
\end{equation*}
\begin{equation*}
v_{k}(x)=\frac{1}{c_{k}}u_{k}(exp_{x_{k}}(r_{k}x));
\end{equation*}
\begin{equation*}
g_{k}(x)=(exp_{x_{k}}^*g)(r_{k}x).
\end{equation*}
Next we define
\begin{equation*}
d_k=\frac{c_k}{\beta_k}\;\;\;\;d=\lim_{k\rightarrow +\infty}d_k.
\end{equation*}
\begin{prop}\label{bubble}
The following hold:\\
We have
\begin{equation*}
\text{if}\;\;\;d<+\infty\;\;\;\text{then}\;\;\;w_k\rightarrow
w(x):=\frac{4}{d}
\log\left(\frac{1}{1+\sqrt{\frac{d}{6}}|x|^2}\right)
\;\;\text{in}\;\; C^2_{loc}(\R^4);
\end{equation*}
and
\begin{equation*}
\text{if}\;\;\;d=\infty\;\;\;\text{then}\;\;\;\;\;\;\;\;w_k\rightarrow
w=0 \;\;\text{in}\;\; C^2_{loc}(\R^4).
\end{equation*}
\begin{proof}
First of all we recall that
\begin{equation*}
g_{k}\rightarrow dx^2\;\;\;\text{in}\;\;C^2_{loc}(\R^4).
\end{equation*}
Since \;$(\frac{\lambda_k}
{\beta_k^2})$,\;$(\frac{\beta_k}{c_k})$\;are bounded and
\;$c_k\rightarrow +\infty$\;, then we infer that
\begin{equation*}
r_k\rightarrow 0\;\;\text{as}\;\;k\rightarrow 0.
\end{equation*}
Now using the Green representation formula for \;$\D_g^2$\;(see
Lemma \;$\ref{eq:greenlap}$) we have that
\begin{equation*}
u_k(x)=\int_{M}F(x,y)\D_g^2 u_kdV_g(y)\;\;\;\forall x\in M.
\end{equation*}
Now using equation and differentiating with respect to \;$x$\; we
obtain that for \;$m=1,2$\;
\begin{equation*}
|\nabla_g^m u_k(x)|\leq \int_{M}\left|\nabla_g^m
F(x,y)\right|\left|\frac{u_k}{\lambda_k}
e^{\alpha_ku_k^2}-\gamma_k\right|dV_g(y).
\end{equation*}
Hence from the fact that \;$\beta_k\gamma_k$\;is bounded see Lemma
\;$\ref{eq:asymli}$\; we get
\begin{equation*}
|\nabla_g^m u_k(x)|\leq \int_{M}\left|\nabla_g^m
F(x,y)\right|\left|\frac{u_k}{\lambda_k}
e^{\alpha_ku_k^2}\right|dV_g(y)+O(\beta_k^{-1}).
\end{equation*}
Now for \;$y_k\in B_{Lr_k}(x_k)$, \;$L>0$\;fixed we write that
\begin{equation*}
\begin{split}
\int_{M}\left|\nabla_g^m
F(y_k,y)\right|\frac{|u_k|}{\lambda_k}e^{\alpha_ku_k^2}dV_g(y)
=O\left(r_k^{-m}\int_{M\setminus
B_{Lr_k}(y_k)}\frac{|u_k|}{\lambda_k}e^{\alpha_ku_k^2}
dV_g(y)\right)\\+O\left(\frac{c_{k}}{\lambda_k}e^{\alpha_kc_k^2}\int_{B_{Lr_k}(y_k)}d_{g}
(y_k,y)^{-m}dV_g(y)\right)\\=O(r_k^{-m}\beta_k^{-1}).
\end{split}
\end{equation*}
thanks to the fact that \;$|u_k|\leq c_{k}$\;to the definition of \;$r_k$.\\
Now it is not worth remarking that \;$c_k=u_k(x_k)$\;since we have
taken \;$\tau\geq 0$\; (see Lemma\;$\ref{eq:asymli}$). Hence we have
that
\begin{equation*}
w_k(x)\leq w_{k}(0)=0\;\;\;\forall \;x\in \R^4.
\end{equation*}
So we get from the estimate above that \;$w_k$\;is uniformly bounded
in \;$C^2(K)$\;for every compact subset \;$K$\;of \;$\R^4$. Thus by
Arz\'ela-Ascoli Theorem we infer that
\begin{equation*}
w_{k}\longrightarrow w \in C^1_{loc}(\R^4).
\end{equation*}
Clearly \;$w$\;is a Lipschitz function since the constant which
bounds the gradient of
\;$w_{k}$\;is independent of the compact set \;$K$.\\
On the other hand from the Green representation formula we have for
\;$x\in \R^4$\; fixed and for \;$L$\;big enough such that \;$x\in
B^L(0)$\;
\begin{equation*}
u_k(exp_{x_{k}}(r_kx))=\int_{M}F(exp_{x_{k}}(r_kx),y)\D_g^2u_k(y)dV_g(y).
\end{equation*}
Now remarking that
\begin{equation*}
u_k(x_{k})=u_k(exp_{x_{k}}(r_k0));
\end{equation*}
we have that
\begin{equation*}
u_k(exp_{x_{k}}(r_kx))-u_k(x_{k})=\int_{M}\left(F(exp_{x_{\epsilon}}(r_kx),y)-F(exp_{x_{k}}
(0),y)\right)\D_g^2u_k(y)dV_g(y).
\end{equation*}
Hence using \;$\eqref{eq:edp}$\;we obtain
\begin{equation*}
\begin{split}
u_k(exp_{x_{k}}(r_kx))-u_k(x_{k})=\int_{M}\left(F(exp_{x_{k}}(r_kx),y)-F(exp_{x_{k}}(0),y)
\right)\frac{u_k}{\lambda_{k}}e^{\alpha_ku_k^2}dV_g(y)\\
-\int_{M}\left(F(exp_{x_{k}}(r_kx),y)
-F(exp_{x_{k}}(0),y)\right)(\gamma_k)dV_g(y).
\end{split}
\end{equation*}
Now setting
\begin{equation*}
I_{k}(x)=\int_{B_{Lr_k}(x_k)}\left(F(exp_{x_{k}}(r_kx),y)-F(exp_{x_{k}}(0),y)\right)\frac{u_k}
{\lambda_{k}}e^{\alpha_ku_k^2}dV_g(y);
\end{equation*}
\begin{equation*}
\text{II}_{k}(x)=\int_{M\setminus
B_{Lr_k}(x_k)}\left(F(exp_{x_{k}}(r_kx),y)-F(exp_{x_{k}}
(0),y)\right)\frac{u_k}{\lambda_{k}}e^{\alpha_ku_k^2}dV_g(y)
\end{equation*}
and
\begin{equation*}
\text{III}_{k}(x)=\int_{M}\left(F(exp_{x_{k}}(r_kx),y)-F(exp_{x_{k}}(0),y)\right)(\gamma_k)
dV_g(y);
\end{equation*}
we find
\begin{equation*}
u_k(exp_{x_{k}}(r_kx))-u_k(x_{k})=I_{k}(x)+\text{II}_{k}(x)+\text{III}_{k}(x).
\end{equation*}
So using the definition of \;$w_{k}$ we arrive to
\begin{equation*}
w_{k}=2\alpha_k\beta_k\left(I_{k}(x)+\text{II}_{k}(x)+\text{III}_{k}(x)\right).
\end{equation*}
Now to continue the proof we consider two cases:\\
{\bf Case 1:\;$d<+\infty$}

\noindent First of all let us study each of the terms
\;$2\alpha_k\beta_kI_{k}(x)$,\;$2\alpha_k
\beta_k\text{II}_{k}(x)$,\;$2\alpha_k\beta_k\text{III}_{k}(x)$\;separately.\\
Using the change of variables \;$y=exp_{x_{k}}(r_kz)$\;we have
\begin{equation*}
\begin{split}
2\alpha_k\beta_kI_{k}(x)=\int_{B^{L}(0)}\left(F(exp_{x_k}(r_kx),exp_{x_{k}}(r_kz))
-F(exp_{x_{k}}(0),exp_{x_{k}}(r_kz))\right)\\\frac{2\alpha_k\beta_ku_k(exp_{x_{k}}(r_kz))}
{\lambda_{k}}e^{\alpha_ku^2_{\epsilon}(exp_{x_{k}}(r_kz)}r_k^4dV_{g_{k}}(z).
\end{split}
\end{equation*}
Hence using the definition of \;$r_k$\; and \;$v_k$\; one can check
easily that the following holds
\begin{equation*}
\begin{split}
2\alpha_k\beta_kI_{k}(x)=2\alpha_k\int_{B^L(0)}\left(G(exp_{x_{\epsilon}}(r_kx),
exp_{x_{\epsilon}}(r_kz))-G(exp_{x_{\epsilon}}(0),exp_{x_{\epsilon}}(r_kz))\right)
v_{k}(z)\\e^{\frac{d_k}{2}(w_{k}(z)(1+v_k)}dV_{g_{k}}(z).
\end{split}
\end{equation*}
Moreover from the asymptotics of the Green function see
Lemma\;$\ref{eq:greenlap}$\; we have that
\begin{equation*}
2\alpha_k\beta_kI_{k}(x)=2\alpha_k\int_{B^L(0)}\left(\frac{1}{8\pi^2}\log\frac{|z|}
{|x-z|}+K_{k}(x,z)\right)v_{k}(z)e^{\frac{d_k}{2}(w_{k}(z)(1+v_k(z)))}dV_{g_{k}}(z).
\end{equation*}
where
\begin{equation*}
K_{k}(x,z)=\left[K(exp_{x_{k}}(r_kx),exp_{x_{k}}(r_kz)\right)-\left(K(exp_{x_{k}}(0),
exp_{x_{k}}(r_kz)\right].
\end{equation*}
Hence since \;$K$\;is of class \;$C^1$\; on \;$M^2$ and
\;$g_{k}\rightarrow dx^2$\;in \;$C^2_{loc}(\R^4)$\;and
\;$v_{k}\rightarrow 1$\;then letting \;$k\rightarrow +\infty$\; we
derive
\begin{equation*}
\lim_{k\rightarrow+\infty}2\alpha_k\beta_kI_{k}(x)=8\int_{B^L(0)}\log\frac{|z|}{|x-z|}e^{dw(z)}dz.
\end{equation*}
Now to estimate\;$\alpha_k\beta_k\text{II}_{k}(x)$\; we write for
\;$k$\;large enough
\begin{equation*}
\begin{split}
\alpha_k\beta_k\text{II}_{k}(x)=\int_{M\setminus
B_{Lr_k}(x_k)}\frac{1}{8\pi^2}\log
\left(\frac{d_{g}(exp_{x_{k}}(0),y)}{d_{g}(exp_{x_{k}}(r_kx),y)}\right)\frac{2\alpha_k
\beta_ku_k}{\lambda_{k}}e^{\alpha_ku_k^2}dV_{g}(y)\\+\int_{M\setminus
B_{Lr_k}(x_k)} \bar
K_{k}(x,y)\frac{2\alpha_k\beta_ku_k}{\lambda_{k}}e^{\alpha_ku_k^2}dV_{g}(y),
\end{split}
\end{equation*}
where
\begin{equation*}
\bar
K_{k}(x,y)=\left(K(exp_{x_{k}}(r_kx),y)-K(exp_{x_{k}}(0),y)\right).
\end{equation*}
Taking the absolute value in both sides of the equality and using
the change of variable \;$y=exp_{x_{k}}(r_kz)$\; and the fact that
\;$K\in C^{1}$\; we obtain,
\begin{equation*}
\begin{split}
|2\alpha_k\beta_k\text{II}_{k}(x)|\leq \int_{\R^4 \setminus
B^{L}(0)}8\left|\log\left(
\frac{|z|}{|x-z|}\right)\right||v_{k|}(z)e^{\frac{d_k}{2}(w_{k}(z)(1+v_k(z)))}dV_{g_{k}}
(z)\\ +Lr_k\int_{M\setminus
B_{Lr_k}(x_k)}\frac{2\alpha_k\beta_ku_k}{\lambda_{k}}
e^{\alpha_ku_k^2}dV_{g}(y).
\end{split}
\end{equation*}
Hence letting \;$k\rightarrow +\infty$\;we deduce that
\begin{equation*}
\limsup_{k\rightarrow+\infty}|2\alpha_k\beta_k\text{II}_{k}(x)|=o_{L}(1).
\end{equation*}
Now using the same method one proves that
\begin{equation*}
2\alpha_k\beta_k\text{III}_{k}(x)\rightarrow 0\;\;as\;\;k\rightarrow
+\infty.
\end{equation*}
So we have that
\begin{equation*}
w(x)=\int_{B^{L}(R)}8
\log\left(\frac{|z|}{|x-z|}\right)e^{dw(z)}dz+\lim_{k\rightarrow+\infty}
2\alpha_k\beta_k\text{II}_{k}(x).
\end{equation*}
Hence letting \;$L\rightarrow+\infty$\;we obtain that \;$w$\; is a
solution of the following integral equation
\begin{equation}\label{eq:inte}
w(x)=\int_{\R^4}8\log\left(\frac{|z|}{|x-z|}\right)e^{dw(z)}dz.
\end{equation}
Now since\;$w$\;is Lipschitz then the theory of singular integral
operator gives
that \;$w\in C^1(\R^4)$.\\
Since
\begin{equation*}
\lim_{k\rightarrow+\infty}\int_{B_{Lr_k}(x_k)}\frac{2\alpha_k\beta_ku_k}{\lambda_{k}}
e^{\alpha_ku_k^2}dV_g
=64\pi^2\int_{B^L(0)}e^{dw(x)}dx.
\end{equation*}
and
\begin{equation*}
\int_{B_{Lr_k}(x_k)}\frac{2\alpha_k\beta_ku_k}{\lambda_{k}}e^{\alpha_ku_k^2}dV_g
\leq 64\pi^2,
\end{equation*}
then we get
\begin{equation*}
\lim_{L\rightarrow+\infty}\int_{B^L(0)}e^{dw(x)}dx=\int_{\R^4}e^{dw(x)}dx\leq
1.
\end{equation*}
Now setting
\begin{equation*}
\tilde w(x)=\frac{d}{4}w(x)+\frac{1}{4}\log(\frac{8\pi^2d}{3});
\end{equation*}
we have that \;$\tilde w$\;satisfies the following conformally
invariant integral equation
\begin{equation}\label{eq:confinte}
\tilde
w(x)=\int_{\R^4}\frac{6}{8\pi^2}\log\left(\frac{|z|}{|x-z|}\right)e^{\tilde
w(z)}dz +\frac{1}{4}\log(\frac{8\pi^2d}{3}),
\end{equation}
and
\begin{equation*}
\int_{\R^4}e^{4\tilde w(x)}dx<+\infty.
\end{equation*}
Hence from the classification result by X.Xu see Theorem 1.2 in
\cite{xu} we derive that
\begin{equation*}
\tilde w(x)=\log\left(\frac{2\lambda}{\lambda^2+|x-x_0|^2}\right)
\end{equation*}
for some \;$\lambda>0$\;and \;$x_0\in \R^4$.\\
>From the fact that
\begin{equation*}
w(x)\leq w(0)=0\;\;\;\forall x\in \R^4;
\end{equation*}
we obtain
\begin{equation*}
\tilde w(x)\leq \tilde
w(0)=\frac{1}{4}\log(\frac{8\pi^2d}{3})\;\;\;\forall x\in \R^4.
\end{equation*}
Then we derive
\begin{equation*}
x_0=0,\;\;\lambda=2(\frac{8\pi^2d}{3})^{-\frac{1}{4}}
\end{equation*}
Hence by trivial calculations we get
\begin{equation*}
w(x)=\frac{4}{d}\log\left(\frac{1}{1+\sqrt{\frac{d}{6}}|x|^2}\right).
\end{equation*}
{\bf Case 2: \;$d=+\infty$.}\\
In this case using the same argument we get
\begin{equation*}
\limsup_{k\rightarrow+\infty}|\alpha_k\beta_k\text{II}_{k}(x)|=o_{L}(1);
\end{equation*}
and
\begin{equation*}
\alpha_k\beta_k\text{III}_{k}(x)=o_{k}(1),
\end{equation*}
Now let us show that
\begin{equation*}
\alpha_k\beta_k\text{I}_{k}(x)=o_{k}(1)
\end{equation*}
By using the same arguments as in Case 1 we get
\begin{equation*}
\alpha_k\beta_kI_{k}(x)=\int_{B^L(0)}\left(\frac{1}{8\pi^2}\log\frac{|z|}{|x-z|}
+K_{k}(x,z)\right)v_{k}(z)e^{d_k(w_{k}(z)(1+v_k(z)))}dV_{g_{k}}(z)
\end{equation*}
Now  since \;$K$\;is \;$C^1$\; we need only to show that
\begin{equation*}
\int_{B^L(0)}\frac{1}{8\pi^2}\log\frac{|z|}{|x-z|}v_{k}(z)e^{d_k(w_{k}(z)(1+v_k(z)))}
dV_{g_{k}}(z)=o_k(1).
\end{equation*}
By using the trivial inequality
\begin{equation*}
\int_{B_{Lr_k}(x_k)}\frac{u_k^2}{\lambda_k}e^{\alpha_ku_k^2}dV_g\leq
1;
\end{equation*}
and the change of variables as above, we obtain
\begin{equation*}
\int_{B^L(0)}v_{k}^2(z)e^{d_k(w_{k}(z)(1+v_k(z)))}dV_{g_{k}}(z)=O(\frac{1}{d_k})=o_{k}(1).
\end{equation*}
On the other hand using the property of \;$v_k$\;one can check
easily that
\begin{equation*}
\int_{B^L(0)}v_{k}(z)e^{d_k(w_{k}(z)(1+v_k(z)))}dV_{g_{k}}(z)=\int_{B^L(0)}v_{k}^2(z)
e^{d_k(w_{k}(z)(1+v_k(z)))}dV_{g_{k}}(z)+o_k(1).
\end{equation*}
Thus  we arrive to
\begin{equation*}
\int_{B^L(0)}\frac{1}{8\pi^2}\log\frac{|z|}{|x-z|}v_{k}(z)e^{d_k(w_{k}(z)(1+v_k(z)))}
dV_{g_{k}}(z)=o_k(1)
\end{equation*}
So we get
\begin{equation*}
\alpha_k\beta_k\text{I}_{k}(x)=o_{k}(1)
\end{equation*}
Thus letting \;$k\rightarrow+\infty$, we obtain
\begin{equation*}
w(x)=0\;\;\forall x\in \R^4.
\end{equation*}
Hence the Proposition is proved.
\end{proof}
\end{prop}

\subsection{Capacity estimates}
This subsection deals with some capacity-type estimates which allow
us to get an upper bound of
\;$\tau^2\lim_{k\rightarrow+\infty}\frac{\lambda_k}{\beta_k^2}$.  We start by giving a first
Lemma to show that we can basically work on Euclidean space in order
to get the capacity estimates as already said in the Introduction.
\begin{lem}\label{eq:capeu}
There is a constant \;$B$\; which is independent
of\;$k$,\;$L$\;and\;$\delta$\; s.t.
\begin{equation*}
\int_{B^{\delta}(0)\setminus
B^{Lr_k}(0)}|(1-B|x|^2)\Delta_0\tilde{u}_k|^2dx \leq
\int_{B_\delta(x_k)\setminus B_{Lr_k}(x_k)}|\Delta_gu_k|^2dV_g+
\frac{J_1(k,L,\delta)}{\beta_k^2},
\end{equation*}
where
\begin{equation*}
\tilde u(x)=u_k(exp_{x_k}(x)).
\end{equation*}
Moreover we have that
\begin{equation*}
\lim_{\delta\rightarrow
0}\lim_{k\rightarrow+\infty}J_1(k,L,\delta)=0.
\end{equation*}
\end{lem}

\begin{proof}
First of all by using the definition of \;$\D_g$\;ie
\begin{equation*}
\D_g=\frac{1}{\sqrt{|g|}}\partial_r(\sqrt{|g|}g^{rs}\partial_s);
\end{equation*}
we get
\begin{equation*}
\begin{array}{lll}
   |\Delta_g\beta_ku_k|^2&=&|g^{rs}\beta_k\frac{\partial^2\tilde{u}_k}
     {\partial x^r\partial x^s}
       +O(|\nabla \beta_k\tilde{u}_k|)|^2\\[\mv]
   &=&|g^{rs}\beta_k\frac{\partial^2\tilde{u}_k}{\partial x^r\partial x^s}|^2+
      O(|\nabla^2\beta_k\tilde{u}_k||\nabla\beta_k \tilde{u}_k|)
       )+O((|\nabla\beta_k \tilde{u}_k|)^2)
\end{array}
\end{equation*}
On the other hand using the fact that (see
Corollary\;$\ref{eq:profile}$))
\begin{equation*}
\beta_k\tilde{u}_k\rightharpoonup \tilde
G\;\;\;\text{in}\;\;W^{2,p}(M);
\end{equation*}
where \;$p\in(1,2)$; and \;$\tilde G(x)=G(exp_{x_0}(x))$; we obtain
\begin{equation*}
\begin{array}{l}
   \dint_{B^\delta(0)\setminus B^{Lr_k}(0)}O(|\nabla^2\beta_k\tilde{u}_k||\nabla
   \beta_k\tilde{u}_k|
    )
    +O((|\nabla\beta_k \tilde{u}_k|)^2)\\[\mv]
     \s\s\s\s\s\s\leq C||\tilde G||_{W^{1,2}(B^\delta(0)\setminus B^{Lr_k}(0))}
      \\[\mv]
     \s\s\s\s\s\s= J_2(k,L,\delta),
  \end{array}
\end{equation*}
and it is clear that
\begin{equation*}
\lim_{\delta\rightarrow
0}\lim_{k\rightarrow+\infty}J_2(k,L,\delta)=0
\end{equation*}
Now let us estimate $ \int_{B^\delta(0)\setminus
B^{Lr_k}(0)}|g^{rs}\beta_k\frac{\partial^2\tilde{u}_k} {\partial
x^r\partial x^s}|^2$.  To do this, we first write the inverse of the
metric in the following form
\begin{equation*}
g^{rs}=\delta^{rs}+A^{rs}
\end{equation*}
with
\begin{equation*}
|A^{rs}|\leq C|x|^2.
\end{equation*}
We can write
\begin{equation*}
|g^{rs}\frac{\partial^2\tilde{u}_k}{\partial x^r\partial
x^s}|^2|\Delta_0\tilde{u}_k|^2+2\sum_{p,q}A^{pq}\D_0\tilde u_k
\frac{\partial^2\tilde{u}_k}{\partial x^p\partial
x^q}+\sum_{r,s,p,q} A^{rs}A^{pq}\frac{\partial^2 \tilde{u}_k}
{\partial
 x^r\partial x^s}
\frac{\partial^2\tilde{u}_k}{\partial x^p\partial x^q}
\end{equation*}
Furthermore we derive
\begin{equation*}
\sum_{p,q}2\int_{B^\delta(0)\setminus B^{Lr_k}(0)}|A^{pq}\D_0\tilde
u_k \frac{\partial^2\tilde{u}_k}{\partial x^p\partial x^q}|dV_g\leq
C\int_{B^\delta(0)\setminus B^{Lr_k}(0)} (|x|^2|\D_0\tilde
u_k|^2+\sum_{p,q}|x|^2| \frac{\partial^2\tilde{u}_k}{\partial
x^p\partial x^q}|^2)dx
\end{equation*}
On the other hand  we have that
\begin{equation*}
\begin{split}
   \sum_{p,q}\int_{B^\delta(0)\setminus B^{Lr_k}(0)}|x|^2|
     \frac{\partial^2\tilde{u}_k}{\partial x^p\partial x^q}|^2dx     \int_{B^\delta(0)\setminus B^{Lr_k}(0)}|x|^2\frac{\partial^2\tilde{u}_k}
      {\partial x^s\partial x^s}
      \frac{\partial^2\tilde{u}_k}{\partial x^p\partial x^p}dx\\
      +\int_{B^\delta(0)\setminus B^{Lr_k}(0)}
     O(|\nabla \tilde{u}_k||\nabla^2 \tilde{u}_k|)dx
   +\int_{\partial(B^\delta(0)\setminus B
   ^{Lr_k}(0))}
     |x|^2\frac{\partial \tilde{u}_k}
      {\partial x^q}
      \frac{\partial^2\tilde{u}_k}{\partial x^p\partial x^q}\left<\frac{\partial}
      {\partial x^p},
      \frac{\partial}{\partial r}\right>dS\\+\int_{\partial(B^\delta(0)\setminus B^{Lr_k}(0))}
      |x|^2\frac{\partial \tilde{u}_k}
      {\partial x^q}
      \frac{\partial^2\tilde{u}_k}{\partial x^p\partial x^p}\left<\frac{\partial}
      {\partial x^q},\frac{\partial}
      {\partial r}\right>dS.
         \end{split}
  \end{equation*}
So setting
\begin{equation*}
\begin{split}
\frac{J_3(k,L,\delta)}{\beta_k^2}=\int_{B^\delta(0)\setminus
B^{Lr_k}(0)}
     O(|\nabla \tilde{u}_k||\nabla^2 \tilde{u}_k|)dx
   +\int_{\partial(B^\delta(0)\setminus B^{Lr_k}(0))}
     |x|^2\frac{\partial \tilde{u}_k}
      {\partial x^q}
      \frac{\partial^2\tilde{u}_k}{\partial x^p\partial x^q}\left<\frac{\partial}
      {\partial x^p},
      \frac{\partial}{\partial r}\right>dS\\+\int_{\partial(B^\delta(0)\setminus B^{Lr_k}(0))}
      |x|^2\frac{\partial \tilde{u}_k}
      {\partial x^q}
      \frac{\partial^2\tilde{u}_k}{\partial x^p\partial x^p})\left<\frac{\partial}
      {\partial x^q},
      \frac{\partial}{\partial r}\right>dS
\end{split}
\end{equation*}
We obtain
\begin{equation*}
\sum_{p,q}\int_{B^\delta(0)\setminus B^{Lr_k}(0)}|x|^2|
     \frac{\partial^2\tilde{u}_k}{\partial x^p\partial x^q}|^2
   =\int_{B^\delta(0)\setminus B^{Lr_k}(0)}|x|^2\frac{\partial^2\tilde{u}_k}
      {\partial x^q\partial x^q}
      \frac{\partial^2\tilde{u}_k}{\partial x^p\partial x^p}dx+\frac{J_3(k,L,\delta)}
      {\beta_k^2}.
\end{equation*}
Moreover we have that
\begin{equation*}
\lim_{\delta\rightarrow
0}\lim_{k\rightarrow+\infty}J_3(k,L,\delta)=0.
\end{equation*}
Hence we get
\begin{equation*}
2\sum_{p,q}\int_{B^\delta(0)\setminus
B^{Lr_k}(0)}|A^{pq}\frac{\partial^2 \tilde{u}_k} {\partial
x^s\partial x^s} \frac{\partial^2\tilde{u}_k}{\partial x^p\partial
x^q}|\leq C\int_{B^\delta(0)\setminus
B^{Lr_k}(0)}|x|^2|\Delta_0\tilde{u}_k|^2dx+\frac{J_4(k,L,\delta)}{\beta_k^2}
\end{equation*}
with
\begin{equation*}
\lim_{\delta\rightarrow
0}\lim_{k\rightarrow+\infty}J_4(k,L,\delta)=0.
\end{equation*}
On the other hand using similar arguments we get
\begin{equation*}\int_{B^\delta(0)\setminus B^{Lr_k}(0)}\sum_{r,s,p,q}
A^{rs}A^{pq}\frac{\partial^2 \tilde{u}_k} {\partial
 x^r\partial x^s}
\frac{\partial^2\tilde{u}_k}{\partial x^p\partial x^q} \leq
C\int_{B^\delta(0)\setminus
B^{Lr_k}(0)}|x|^4|\Delta_0\tilde{u}_k|^2dx+\frac{J_5(k,L,\delta)}{\beta_k^2}.
\end{equation*}
with
\begin{equation*}
\lim_{\delta\rightarrow
0}\lim_{k\rightarrow+\infty}J_5(k,L,\delta)=0.
\end{equation*}
So we arrive to
\begin{equation*}
\int_{B_\delta(x_k)\setminus B_{Lr_k}(x_k)}|\Delta_gu_k|^2dV_g
     \leq\int_{B^\delta(0)\setminus B^{Lr_k}(0)}(1+C|x|^2+C|x|^4)|\Delta_0\tilde{u}_k|^2dx
     +\frac{J_6(k,L,\delta)}{\beta_k^2};
\end{equation*}
with
\begin{equation*}
\lim_{\delta\rightarrow
0}\lim_{k\rightarrow+\infty}J_6(k,L,\delta)=0
\end{equation*}
Hence we can find a constant \;$B_1$\;independent of
\;$k$,\;$L$\;and\;$\delta$ s.t
\begin{equation*}
\int_{B_\delta(x_k)\setminus B_{Lr_k}}|\Delta_gu_k|^2dV_g
     \geq\int_{B^\delta(0)\setminus B^{Lr_k}(0)}(1-B_1|x|^2)|\Delta_0\tilde{u}_k|^2dx+
     \frac{J_7(k,L,\delta)}{\beta_k^2}.
\end{equation*}
So setting
\begin{equation*}
   J_1(k,L,\delta)=-J_7(k,L,\delta)\;\;\;\text{and}\;\;B=B_1
\end{equation*}
we have the proved the Lemma.
\end{proof}
Next we give a technical Lemma
\begin{lem}\label{eq:}
There exists a sequence of functions \; $U_k\in
W^{2,2}(B^\delta(0)\setminus B^{Lr_k}(0))$\; s.t
\begin{equation*}
U_k|_{\partial
B^\delta(0)}=\tau\frac{-\frac{1}{16\pi^2}\log\delta+S_0
}{\beta_k},\s U_k|_{\partial
B^{Lr_k}(0)}=\frac{w(L)}{2\alpha_k\beta_k}+c_k;
\end{equation*}
and
\begin{equation*}
\frac{\partial U_k}{\partial r}|_{
\partial B_\delta(0)}=-\frac{\tau}{8\pi^2\delta\beta_k},\s
\frac{\partial U_k}{\partial r}|_{\partial
B^{Lr_k}(0)}=\frac{w'(L)}{2\alpha_k\beta_kr_k}.
\end{equation*}
Moreover there holds
\begin{equation*}
\lim_{\delta\rightarrow 0}\lim_{k\rightarrow+\infty}\beta_k^2
(\int_{B^\delta(0)\setminus B^{Lr_k}(0)}|\Delta_0(1-B|x|^2)U_k|^2dx-
\int_{B^\delta\setminus
B^{Lr_k}(0)}|(1-B|x|^2)\Delta_0\tilde{u}_k|^2dx)=0.
\end{equation*}
\end{lem}
\begin{proof}
First of all let us set
\begin{equation*}
h_k(x)=u_k(exp_{x_k}(r_kx)).
\end{equation*}
and \; $u_k'$\; to  be the solution of
\begin{equation*}
\left\{\begin{array}{l}
             \Delta_0^2u_k'=\Delta_0^2h_k\\[\mv]
               \frac{\partial u_k'}{\partial n}|_{\partial B_{2L}}                \frac{\partial h_k}{\partial n}|_{\partial B_{2L}},\s
                 u_k'|_{\partial B^{2L}(0)}=h_k|_{\partial B^{2L}(0)}\\[\mv]
               \frac{\partial u_k'}{\partial n}|_{\partial B^{L}(0)}                \frac{1}{2\alpha_k\beta_k}\frac{\partial w}{\partial n}|_{\partial B^{L}(0)},\s
                 u_k'|_{\partial B^{L}(0)}=\frac{w}{2\alpha_k\beta_k}|_{\partial B^{L}(0)}.
          \end{array}\right.
\end{equation*}
Next let us define
\begin{equation*}
U_k'=\left\{\begin{array}{ll}
                        u_k'(\frac{x}{r_k})&Lr_k\leq
                          |x|
                          \leq 2Lr_k\\[\mv]
                        \tilde{u}_k(x)&2Lr_k\leq
                          |x|.
                     \end{array}\right.
\end{equation*}
Clearly we have that
\begin{equation*}
\lim_{k\rightarrow +\infty}\int_{B^{2Lr_k}(0)\setminus B^{Lr_k}(0)}
(1-B|x|^2)(|\Delta_0U_k'|^2-|\Delta_0\tilde{u}_k|^2)dx=0,$$ and
$$\lim_{k\rightarrow+\infty}
|U_k'-\tilde{u}_k'|_{C^0(B^{2Lr_k}(0)\setminus B^{Lr_k}(0))}=0.
\end{equation*}
Now let \;$\eta$\;be a smooth function which satisfies
\begin{equation*}
\eta(t)=\left\{\begin{array}{ll}
                      1&t\leq 1/2\\[\mv]
                      0& t>2/3
                 \end{array}\right.
\end{equation*}
and set
\begin{equation*}
  G_k=\eta(\frac{|x|}{\delta})(\tilde{u}_k-\tau S_0
+\frac{\tau}{8\pi^2}\log{|x|}) -\frac{\tau}{8\pi^2}\log{|x|}+\tau
S_0.
\end{equation*}
Then we have that
\begin{equation*}
G_k\rightarrow -\frac{\tau}{8\pi^2}\log{|x|}+\tau
S_0+\tau\eta(\frac{|x|}{\delta})\tilde S_1(x);
\end{equation*}
where \;$\tilde S_1(x)=S_1(exp_{x_0}(x))$\;.\\
Furthermore we obtain
\begin{equation*}
 \beta_k\tilde{u}_k-G_k
\rightarrow \tau\left(1-\eta(\frac{|x|}{\delta})\right)S_1(x),
\end{equation*}
then
\begin{equation*}
\lim\limits_{\epsilon\rightarrow 0}
  |\int_{B^\delta(0)\setminus B^{\delta/2}(0)}|\Delta_0\beta_k\tilde{u}_k|^2dx-
    \int_{B^\delta(0)\setminus B^{\delta/2}(0)}|\Delta_0G_k|^2dx|
    \leq \Sigma.
\end{equation*}
    where
\begin{equation*}
\begin{array}{lll}
    \Sigma&=&\sqrt{\int_{B^\delta(0)\setminus B^{\delta/2}(0)}|\Delta_0(1-\eta(\frac{|x|}
       {\delta}))\tilde S_1(x)|^2dx
       \int_{B^\delta(0)\setminus B^{\delta/2}(0)}|\Delta_0(\tilde G-\frac{1}{8\pi^2}
     \log{|x|}+\eta(\frac{|x|}{\delta})\tilde S_1(x))|^2dx}\\[\mv]
   &\leq& C\delta\sqrt{|\log\delta|}.
\end{array}
\end{equation*}
So we get
\begin{equation*}
\lim\limits_{\epsilon\rightarrow 0}
    |\int_{B^\delta(0)\setminus B^{\delta/2}(0)}|\Delta_0\beta_k\tilde{u}_k|^2dx-
    \int_{B^\delta(0)\setminus B^{\delta/2}(0)}|\Delta_0G_k|^2dx|
    \leq C\delta\sqrt{|\log\delta|}.
\end{equation*}
Hence setting
\begin{equation*}
U_k=\left\{\begin{array}{ll}
                        U_k'(x)&|x|\leq\frac{\delta}{2}\\[\mv]
                        G_k(x)&\delta/2\leq
                          |x|
                          \leq \delta
                     \end{array}\right.
             \end{equation*}
we have  proved the Lemma.
\end{proof}
\begin{prop}\label{eq:upperbound}
We have the following holds
\begin{equation*}
\tau^2\lim_{k\rightarrow+\infty}\frac{\lambda_k}{\beta_k^2} \leq
\frac{\pi^2}{6}e^{\frac{5}{3}+32\pi^2S_0};
\end{equation*}
and
\begin{equation*}
d\tau=1.
\end{equation*}
\end{prop}
\begin{proof}
First using Lemma \ref{eq:capeu} and Lemma \ref{eq:} we get
\begin{equation}\label{cap1}
\int_{B^\delta(0)\setminus B^{Lr_k}(0)}|\Delta_0(1-B|x|^2)
U_k|^2dx\leq 1-\frac{\int_{B_L(x_0)}|\Delta w|^2+\int_{M\setminus
B_\delta(x_0)}|\Delta G|^2 +J_0(k,L,\delta)} {\beta_k^2}.
\end{equation}
with
\begin{equation*}
\lim_{\delta\rightarrow
0}\lim_{k\rightarrow+\infty}J_0(k,L,\delta)=0.
\end{equation*}
Next we will apply capacity to give a lower boundary of\;
$\int_{B^\delta(0)\setminus B^{Lr_k}(0)}
|\Delta_0(1-B|x|^2)U_k|^2dx$.\\
Hence we need to calculate
$$\inf_{\Phi|_{\partial B^r(0)}=P_1,\Phi|_{\partial B^R(0)}=P_2,
\frac{\partial\Phi}{\partial r}|_{\partial B^r(0)}=Q_1,
\frac{\partial\Phi}{\partial r}|_{\partial B^R(0)}=Q_2}
\int_{B^R(0)\setminus B^r(0)}|\Delta_0\Phi|^2dx,$$
where $P_1$, $P_2$, $Q_1$, $Q_2$ are constants. \\
It is obvious that the infimum is attained by the function
\;$\Phi$\; which satisfies
\begin{equation*}
\left\{\begin{array}{l}
             \Delta^2_0\Phi=0\\[\mv]
             \Phi|_{\partial B^r(0)}=P_1\,\,\,,\Phi|_{\partial B^R(0)}=P_2\,\,\,,
               \frac{\partial\Phi}{\partial r}|_{\partial B^r(0)}=Q_1\,\,\,,
                \frac{\partial\Phi}{\partial r}|_{\partial B^R(0)}=Q_2.
         \end{array}\right.
\end{equation*}
Moreover we can require the function\;$\Phi$\;to be of the form
\begin{equation*}
\Phi=A\log{r}+Br^2+\frac{C}{r^2}+D,
\end{equation*}
where \;$A$,\;$B$, \;$C$, \;$D$\; are all constants which satisfies
the following linear system of equations
\begin{equation*}
\left\{\begin{array}{l}
              A\log{r}+Br^2+\frac{C}{r^2}+D=P_1\\[\mv]
              A\log{R}+BR^2+\frac{C}{R^2}+D=P_2\\[\mv]
              \frac{A}{r}+2Br-2\frac{C}{r^3}=Q_1\\[\mv]
              \frac{A}{R}+2BR-2\frac{C}{R^3}=Q_2
          \end{array}\right.
\end{equation*}
Now by straightforward calculations we obtain the explicit
expression of \;$A$\;and \;$B$\;
\begin{equation*}
\left\{\begin{array}{l}
            A=\frac{P_1-P_2+\frac{\varrho}{2} rQ_1+\frac{\varrho}{2} RQ_2}
             {\log{r/R}+\varrho}\\[\mv]
            B=\frac{-2P_1+2P_2-rQ_1(1+\frac{2r^2}{R^2-r^2}\log{r/R})
                +RQ_2(1+\frac{2R^2}{R^2-r^2}\log{r/R})}{4(R^2+r^2)(\log{r/R}+\varrho)}
          \end{array}\right.
\end{equation*}
Where $\varrho=\frac{R^2-r^2}{R^2+r^2}$. Furthermore we have
\begin{equation}\label{cap2}
\int_{B^R(0)\setminus B^r(0)}|\Delta_0\Phi|^2dx=-8\pi^2
A^2\log{r/R}+
       32\pi^2 AB(R^2-r^2)+32\pi^2B^2(R^4-r^4)
\end{equation}

In our case in which we have that
\begin{equation*}
R=\delta \;\;\;\;r=Lr_k,
\end{equation*}
\begin{equation*}
P_1=c_k+\frac{w(L)}{2\alpha_k\beta_k}+O(r_kc_k)\;\;\;P_2\frac{-\frac{\tau}{8\pi^2}\log\delta+\tau
S_0+O(\delta\log\delta)}{\beta_k}
\end{equation*}
\begin{equation*}
Q_1=\frac{w'(L)+O(r_kc_k)}{2\alpha_k\beta_kr_k}\;\;\;
Q_2=-\frac{\tau+O(\delta\log\delta)}{8\pi^2\beta_k\delta}.
\end{equation*}
Then by the formula giving \;$A$\;we obtain by trivial calculations
\begin{equation*}
A=\frac{c_k+\frac{N_k+\frac{\tau}{8\pi^2}\log\delta}{\beta_k}}{-\log\delta+\log
L+ \frac{\log\frac{\lambda_k}{\beta_kc_k}
-\alpha_kc_k^2}{4}+1+O(r_k^2)}
\end{equation*}
where
\begin{equation*}
N_k=\frac{w(L)}{2\alpha_k}-\tau S_0+ \frac{
w'(L)L}{4\alpha_k}-\frac{\tau}{16\pi^2}+O(\delta\log\delta)+O(r_kc_k^2).
\end{equation*}
Moreover using the the fact that the
sequence\;$(\frac{\lambda_k}{\beta_k^2})_k$\;is bounded it is easily
seen that
\begin{equation*}
  A=O(\frac{1}{c_k}).
  \end{equation*}
Furthermore using the formula of \;$B$\;we get still by trivial
calculations
\begin{equation*}
    B=\frac{-2c_k+\frac{\alpha_kc_k^2}{8\pi^2\beta_k}\frac{\tau}{2}+O(\frac{1}{\beta_k})}
          {\delta^2(-\alpha_kc_k^2+\log\frac{\lambda_k}{\beta_kc_k})}.
\end{equation*}
and then
\begin{equation*}
     B=O(\frac{1}{\beta_k})\frac{1}{\delta^2}.
\end{equation*}
Now let compute\;$ 8\pi^2A^2\log r/R$. By using the expression of
\;$A$,\;$r$\; and \;$R$\;, we have that
\begin{equation*}
  -8\pi^2A^2\log( \frac{r}{R})=-8\pi^2(\frac{c_k+\frac{N_k+\frac{\tau}{8\pi^2}\log\delta}
  {\beta_k}}{-\log\delta+\log L+
     \frac{\log\frac{\lambda_k}{\beta_kc_k}
         -\alpha_kc_k^2}{4}+1+O(r_k^2)})^2(\frac{\log\frac{\lambda_k}{\beta_kc_k}
         -\alpha_kc_k^2}{4}-\log\delta+\log L)
     \end{equation*}
Now using the relation
\begin{equation*}
\begin{split}
(\frac{\alpha_kc_k^2}{4})^2\left(1-\frac{1}{\alpha_kc_k^2}(-4\log\delta+4\log
L+
     \log\frac{\lambda_k}{\beta_kc_k}
         +4+O(r_k^2))\right)^{2}=\\\left(-\log\delta+\log L+
     \frac{\log\frac{\lambda_k}{\beta_kc_k}
         -\alpha_kc_k^2}{4}+1+O(r_k^2)\right)^2
     \end{split}
\end{equation*}
we derive
\begin{equation*}
\begin{split}
 -8\pi^2A^2\log( \frac{r}{R})= -8\pi^2(\frac{c_k+\frac{N_k+\frac{\tau}{8\pi^2}\log\delta}
 {\beta_k}}{
     \frac{\alpha_kc_k^2}{4}})^2\left(1-\frac{1}{\alpha_kc_k^2}(-4\log\delta+4\log L+
     \log\frac{\lambda_k}{\beta_kc_k}
         +4+O(r_k^2))\right)^{-2}\\
         \times(\frac{\log\frac{\lambda_k}{\beta_kc_k}
         -\alpha_kc_k^2}{4}-\log\delta+\log L).
     \end{split}
     \end{equation*}
On the other hand using Taylor expansion we have the following
identity
     \begin{equation*}
     \begin{split}
     \left(1-\frac{1}{\alpha_kc_k^2}(-4\log\delta+4\log L+
     \log\frac{\lambda_k}{\beta_kc_k}
         +4+O(r_k^2))\right)^{-2}=1+
         2\frac{\log\frac{\lambda_k}{\beta_kc_k}+4-4\log\delta+4\log L}{\alpha_kc_k^2}\\
          +O(\frac{\log^2 c_k}{c_k^4});
      \end{split}
     \end{equation*}
     hence we get
     \begin{equation*}
     \begin{split}
 -8\pi^2A^2\log( \frac{r}{R}) =-8\pi^2(\frac{c_k+\frac{N_k+\frac{\tau}{8\pi^2}\log\delta}
 {\beta_k}}{
         \frac{\alpha_kc_k^2}{4}})^2(\frac{\log\frac{\lambda_k}{\beta_kc_k}
         -\alpha_kc_k^2}{4}-\log\delta+\log L)\\
  \times(1+
         2\frac{\log\frac{\lambda_k}{\beta_kc_k}+4-4\log\delta+4\log L}{\alpha_kc_k^2}
          +O(\frac{\log^2 c_k}{c_k^4}))
     \end{split}
     \end{equation*}

On the other hand using the relation
     \begin{equation*}
     \begin{split}
     -8\pi^2(\frac{c_k+\frac{N_k+\frac{\tau}{8\pi^2}\log\delta}{\beta_k}}{
         \frac{\alpha_kc_k^2}{4}})^2(\frac{\log\frac{\lambda_k}{\beta_kc_k}
         -\alpha_kc_k^2}{4}-\log\delta+\log L)=\\\frac{32\pi^2}{\alpha_k}\frac{1}{c_k^2}
         (c_k+\frac{N_k+\frac{\tau}{8\pi^2}\log\delta}
  {\beta_k})^2(1-\frac{\log\frac{\lambda_k}{\beta_kc_k}
         -4\log\delta+4\log L}{\alpha_kc_k^2})
     \end{split}
     \end{equation*}
     we obtain
     \begin{equation*}
     \begin{split}
 -8\pi^2A^2\log( \frac{r}{R}) =\frac{32\pi^2}{\alpha_k}\frac{1}{c_k^2}(c_k+\frac{N_k
 +\frac{\tau}{8\pi^2}\log\delta}
  {\beta_k})^2(1+
         2\frac{\log\frac{\lambda_k}{\beta_kc_k}+4-4\log\delta+4\log L}{\alpha_kc_k^2}
          +O(\frac{\log^2 c_k}{c_k^4}))\\
  \times(1-\frac{\log\frac{\lambda_k}{\beta_kc_k}
         -4\log\delta+4\log L}{\alpha_kc_k^2})
     \end{split}
     \end{equation*}
     Moreover using again the trivial relation
     \begin{equation*}
     \begin{split}
     (1+
         2\frac{\log\frac{\lambda_k}{\beta_kc_k}+4-4\log\delta+4\log L}{\alpha_kc_k^2}
          +O(\frac{\log^2 c_k}{c_k^4}))(1-\frac{\log\frac{\lambda_k}{\beta_kc_k}
         -4\log\delta+4\log L}{\alpha_kc_k^2})=\\(1+
         \frac{\log\frac{\lambda_k}{\beta_kc_k}+8-4\log\delta+4\log L}{\alpha_kc_k^2}
          +O(\frac{\log^2 c_k}{c_k^4}))
     \end{split}
     \end{equation*}
     we arrive to
     \begin{equation*}
  -8\pi^2A^2\log( \frac{r}{R})=\frac{32\pi^2}{\alpha_k}\frac{1}{c_k^2}(c_k+\frac{N_k
  +\frac{\tau}{8\pi^2}\log\delta}{\beta_k})^2
       (1+
         \frac{\log\frac{\lambda_k}{\beta_kc_k}+8-4\log\delta+4\log L}{\alpha_kc_k^2}
          +O(\frac{\log^2 c_k}{c_k^4}))
      \end{equation*}
      On the other hand one can check easily that the following holds
      \begin{equation*}
      \begin{split}
      (c_k+\frac{N_k+\frac{\tau}{8\pi^2}\log\delta}{\beta_k})^2
       (1+
         \frac{\log\frac{\lambda_k}{\beta_kc_k}+8-4\log\delta+4\log L}{\alpha_kc_k^2}
          +O(\frac{\log^2 c_k}{c_k^4}))=\\\left(c_k^2+\frac{\log\frac{\lambda_k}{\beta_kc_k}
          +8-4\log\delta+4\log L}{\alpha_k}
   +2c_k\frac{N_k+\frac{\tau}{8\pi^2}\log\delta}{\beta_k}+O(\frac{\log c_k}{c_k^2})+O(\frac{1}
   {\beta_k^2})\right);
      \end{split}
      \end{equation*}
      thus we obtain
      \begin{equation*}
      \begin{split}
 -8\pi^2A^2\log( \frac{r}{R}) =\frac{32\pi^2}{\alpha_k}\frac{1}{c_k^2}
   \left(c_k^2+\frac{\log\frac{\lambda_k}{\beta_kc_k}+8-4\log\delta+4\log L}{\alpha_k}
   +2c_k\frac{N_k+\frac{\tau}{8\pi^2}\log\delta}{\beta_k}\right)\\+\frac{32\pi^2}{\alpha_k}
   \frac{1}{c_k^2}
   \left(O(\frac{\log c_k}{c_k^2})+O(\frac{1}{\beta_k^2})\right)
   \end{split}
   \end{equation*}
   Furthermore using the relation
   \begin{equation*}
   \begin{split}
   \left(c_k^2+\frac{\log\frac{\lambda_k}{\beta_kc_k}+8-4\log\delta+4\log L}{\alpha_k}
   +2c_k\frac{N_k+\frac{\tau}{8\pi^2}\log\delta}{\beta_k}+O(\frac{\log c_k}{c_k^2})
   +O(\frac{1}{\beta_k^2})\right)=\\
   \left(c_k^2+
      \frac{1}{\alpha_k}\log\frac{\lambda_k}{\beta_kc_k}
      -\frac{4}{\alpha_k}\log\delta+\frac{1}{4\pi^2}d_k\tau\log\delta
      +2d_kN_k+\frac{4\log L}{\alpha_k}
      +\frac{8}{\alpha_k}+o_k(1)\right)
   \end{split}
   \end{equation*}
   we get
   \begin{equation}\label{cap3}
   \begin{split}
 -8\pi^2A^2\log( \frac{r}{R})=\frac{32\pi^2}{\alpha_k^2}\frac{1}{c_k^2}
     \left(c_k^2+
      \frac{1}{\alpha_k}\log\frac{\lambda_k}{\beta_kc_k}
      -\frac{4}{\alpha_k}\log\delta+\frac{1}{4\pi^2}d_k\tau\log\delta
      +2d_kN_k+\frac{4\log L}{\alpha_k}
      +\frac{8}{\alpha_k}\right)\\+\frac{32\pi^2}{\alpha_k^2}\frac{1}{c_k^2}o_k(1)
      \end{split}
  \end{equation}
Next we will evaluate\;$\int_{M\setminus
B_\delta(x_0)}\Delta_gG\Delta_gGdV_g$. We have that by Green formula
\begin{equation*}
    \int_{M\setminus B_\delta(x_0)}\Delta_gG\Delta_gGdV_g=\int_{M\setminus B_\delta(x_0)}G\Delta_g^2G dV_g-\int_{\partial B_\delta(x_0)}
       \frac{\partial G}{\partial r}\Delta_g G+\int_{\partial B_\delta(x_0)}
       G\frac{\partial \Delta_g G}{\partial r}.
       \end{equation*}
       Thus using the equation solved by \;$G$\;we get
       \begin{equation*}
       \begin{split}
  \int_{M\setminus B_\delta(x_0)}\Delta_gG\Delta_gGdV_g  =-\frac{\tau}{\mu(M)}
  \int_{M\setminus B_\delta(p)}GdV_g-
      \frac{\tau^2}{64\pi^4}\int_{\partial B_\delta(x_0)}\frac{\partial(-\log r)}{\partial r}
       \Delta_0(-\log r)\\
    +\int_{\partial B_\delta(x_0)}(-\frac{\tau}{8\pi^2}\log r+S_0)\frac{\partial\Delta_0
     (-\frac{\tau}{8\pi^2}\log r)}{\partial r}
    +O(\delta\log\delta)
    \end{split}
    \end{equation*}
    Hence we obtain
    \begin{equation*}
  \int_{M\setminus B_\delta(x_0)}\Delta_gG\Delta_gGdV_g  =-\frac{\tau^2}{16\pi^2}-
  \frac{\tau^2}{8\pi^2}\log\delta+\tau^2S_0+O(\delta\log\delta),
  \end{equation*}
Now let us set
\begin{equation*}
 P(L)=\int_{B^L(0)}|\Delta_0w|^2dx/(2\times 32\pi^2)^2.
\end{equation*}
Hence using (\ref{cap1}), (\ref{cap2}), (\ref{cap3}), we derive that
\begin{equation*}
\begin{split}
     \frac{32\pi^2}{\alpha_k}\left(c_k^2+
      \frac{1}{\alpha_k}\log\frac{\lambda_k}{\beta_kc_k}
      -\frac{4}{\alpha_k}\log\delta+\frac{1}{4\pi^2}d_k\tau\log\delta
      +2d_kN_k+\frac{4\log L}{\alpha_k}
      +\frac{8}{\alpha_k}\right)\\
   \leq c_k^2(1-\frac{P(L)-\frac{\tau^2}{16\pi^2}
      -\frac{\tau^2}{8\pi^2}\log\delta+\tau S_0+O(\delta\log\delta)+o_{k,\delta}(1)
        }{\beta_k^2})+\delta^2O(c_k^2AB)+\delta^4O(c_k^2B^2).
\end{split}
\end{equation*}
Moreover by isolating the term
\;$\frac{32\pi^2}{\alpha_k^2}\log\frac{\lambda_k}{\beta_kc_k}$\; in
the left and transposing all the other in the right we get
\begin{equation}\label{llbc}
\begin{split}
     \frac{32\pi^2}{\alpha_k^2}\log\frac{\lambda_k}{\beta_kc_k}
     \leq \frac{1}{8\pi^2}(d_k^2\tau^2-\frac{64}{\alpha_k}d_k\tau
     +(\frac{32\pi}{\alpha_k})^2)\log\delta-\frac{32\pi^2}{\alpha_k}
      (2d_kN_k+\frac{4\log L}{\alpha_k}
      +\frac{8}{\alpha_k})\\
     -d_k^2(P(L)+\tau S_0-\frac{\tau^2}{16\pi^2}+O(\delta\log\delta)+o_k(1))
     +\delta^2O(c_k^2AB)+\delta^4O(c_k^2B^2).
  \end{split}
  \end{equation}
  Hence using the trivial identity
  \begin{equation*}
  \log\frac{\lambda_k}{\beta_k^2}=\log\frac{\lambda_k}{\beta_kc_k}+\log d_k
  \end{equation*}
  we get
  \begin{equation*}
\begin{split}
     \frac{32\pi^2}{\alpha_k^2}\log\frac{\lambda_k}{\beta_k^2}
     \leq \frac{1}{8\pi^2}(d_k^2\tau^2-\frac{64}{\alpha_k}d_k\tau
     +(\frac{32\pi}{\alpha_k})^2)\log\delta-\frac{32\pi^2}{\alpha_k}
      (2d_kN_k+\frac{2+4\log L}{\alpha_k}
      +\frac{2}{\alpha_k})\\
     -d_k^2(P(L)+\tau S_0-\frac{\tau^2}{16\pi^2}+O(\delta\log\delta)+o_k(1))+\frac{32\pi^2}
     {\alpha_k^2}\log d_k+O(d_k^2).
  \end{split}
  \end{equation*}
  Now suppose \;$d=+\infty$, letting $\delta\rightarrow 0$, then we have that
  \begin{equation*}
  \lim_{k\rightarrow+\infty} \log\frac{\lambda_k}{\beta_k^2}=-\infty,
  \end{equation*}
  thus we derive
  \begin{equation*}
  \lim_{k\rightarrow+\infty} \frac{\lambda_k}{\beta_k^2}=0
  \end{equation*}
  Hence using Corollary \ref{eq:sup}  we obtain a contradiction. So \;$d$\;must be finite.\\
On the other hand one can check easily that the following holds
\begin{equation*}
\frac{32\pi^2}{\alpha_k^2}\log\frac{\lambda_k}{\beta_k^2} \leq
\frac{1}{8\pi^2}(d_k\tau-\frac{32\pi^2}{\alpha_k})^2\log\delta+O(1)(d_k^2+d_k+\log
d_k)+O(1).
\end{equation*}
Hence we derive
\begin{equation*}
d_k\tau\rightarrow 1;
\end{equation*}
otherwise we reach the same contradiction. So we have that
\begin{equation*}
d\tau=1.
\end{equation*}
Hence by using this we can rewrite \;$B$\;as follows
\begin{equation*}
B=\frac{-2c_k+\delta(-\frac{1}{8\pi^2c_k\delta}2\frac{-\alpha_kc_k^2}{4})+O(1/c_k)}
{\delta^2(-\alpha_kc_k^2)+O(1)}=\frac{o_k(1)}{c_k}.
\end{equation*}
Thus we obtain
\begin{equation*}
32\pi^2AB(R^2-r^2)+32\pi^2B^2(R^4-r^4)=\frac{o_k(1)}{c_k^2}.
\end{equation*}
On the other hand since \;$d<+\infty$, we have that by Lemma
\ref{bubble}
\begin{equation*}
w=-\frac{4\log(1+\sqrt{\frac{d}{6}}\pi|x|^2)}{d}.
\end{equation*}
Moreover by trivial calculations we get
\begin{equation*}
P(L)=\frac{1}{96d^2\pi^2} +\frac{\log(1+\sqrt{\frac{d}{6}}\pi
L^2)}{16d^2\pi^2}.
\end{equation*}
Furthermore by taking the limit as \;$k\rightarrow+\infty$\;in
(\ref{llbc}) we obtain
\begin{equation*}
\lim_{k\rightarrow+\infty}\log\frac{\lambda_k}{\beta_kc_k} \leq
-\frac{25}{3}+4d\tau+2d^2\tau^2+32\pi^2S_0+\frac{4\sqrt{\frac{d}{6}}\pi
L^2}{1+\sqrt{\frac{d}{6}}\pi L^2} +2\log(1+\sqrt{\frac{d}{6}}\pi
L^2)-4\log L
\end{equation*}
Now letting\;$L\rightarrow+\infty$, we get
\begin{equation*}
\lim_{k\rightarrow+\infty}\log\frac{\lambda_k}{\beta_kc_k} \leq
\frac{5}{3}-\log 6 +\log\pi^2 +\log d.
\end{equation*}
Hence by remarking the trivial identity
\begin{equation*}
\lim_{k\rightarrow+\infty}\frac{\lambda_k}{\beta_kc_k}\frac{1}{d}\lim_{k\rightarrow+\infty}\frac{\lambda_k}{\beta_k^2}
\end{equation*}
we get
\begin{equation*}
\tau^2\lim_{k\rightarrow+\infty}\frac{\lambda_k}{\beta_k^2} \leq
\frac{\pi^2}{6}e^{\frac{5}{3}+32\pi^2S_0}.
\end{equation*}
So the proof of the proposition is done.
\end{proof}
\subsection{The test function}
This Subsection deals with the construction of some test functions
in order to reach a
contradiction. \\
Now let \;$\epsilon>0$,\;$c>0$,\;$L>0$\;and set
$$f_\epsilon(x)=\left\{
      \begin{array}{ll}
           c+\frac{\Lambda+Bd_g(x,x_0)^2-4\log\left(1+\lambda (\frac{d_g(x,x_0)}{\epsilon})^2
           \right)}{64\pi^2c}+\frac{S(x)}{c}&d_g(x,x_0)\leq L\epsilon\\[\mv]
          \frac{G(x)}{c}&d_g(x,x_0)>L\epsilon
     \end{array}\right.$$
where
\begin{equation*}
\lambda=\frac{\pi}{\sqrt{6}},\s B=-\frac{4}{L^2\epsilon^2(1+\lambda
L^2)}
\end{equation*}
and
\begin{equation}\label{Lambda}
\Lambda=-64\pi^2c^2-BL^2\epsilon^2-8\log(L\epsilon)+4\log(1+\lambda
L^2).
\end{equation}
\begin{prop}\label{eq:test}
We have that for \;$\epsilon$\;small, there exist suitable
\;$c$\;and \;$L$\;such that
\begin{equation*}
\int_M|\Delta_gf_\epsilon|^2dV_g=1;
\end{equation*}
and
\begin{equation*}
\limsup_{\epsilon\rightarrow
0}\int_Me^{32\pi^2(f_\epsilon-\bar{f}_\epsilon)^2}dV_g
   >Vol(M)+\frac{\pi^2}{6}e^{\frac{5}{3}+32\pi^2S_0}.
\end{equation*}
\end{prop}
\begin{proof}
First of all using the expansion of \;$g$\;in normal coordinates we
get
\begin{equation*}
     \int_{B_{L\epsilon}(x_0)}|\Delta_gf_\epsilon|^2dV_g         \int_{B^{L\epsilon}(0)}|\Delta_0\tilde f_\epsilon|^2(1+O(L\epsilon)^2)dx
        + \int_{B^{L\epsilon}(0)}O(r^2|\nabla_0\tilde f_\epsilon|^2)dx
\end{equation*}
where
\begin{equation*}
\tilde f_\epsilon(x)=f_\epsilon(exp_{x_0}(x)).
\end{equation*}
On the other hand by direct calculations owe obtain
\begin{equation*}
\int_{B^{L\epsilon}(0)}|\Delta_0\tilde
f_\epsilon|^2dx=\frac{12+\lambda L^2 (30+\lambda L^2(21+\lambda
L^2))+6(1+\lambda L^3)^3\log(1+\lambda L^2)}
        {96c^2(1+\lambda L^2)^3\pi^2}
\end{equation*}
Hence we arrive to
\begin{equation*}
\begin{array}{lll}
  \int_{B_{L\epsilon}(x_0)}|\Delta_gf_\epsilon|^2dV_g &=&(1+O(L\epsilon)^2)
  \frac{12+\lambda L^2(30+\lambda L^2(21+\lambda L^2))+6(1+\lambda L^3)^3\log(1+\lambda L^2)}
        {96c^2(1+\lambda L^2)^3\pi^2}\\[\mv]
  &=&\frac{\frac{1}{3}
  +4\log(1+\lambda L^2)+O(\frac{1}{L^2})+O((L\epsilon)^2\log L\epsilon)}{32c^2\pi^2}
\end{array}
\end{equation*}
Furthermore, by direct computation, we have
\begin{equation*}
 \int_{B^{L\epsilon}(0)}r^2|\nabla_0\tilde f_\epsilon|^2dx=O(\frac{L^4\epsilon^4}{c^2}).
\end{equation*}
Moreover using Green formula we get
$$\begin{array}{lll}
  \dint_{M\setminus B_{L\epsilon}(x_0)}|\Delta_gG|^2dV_g
    &=&\dint_{M\setminus B_{L\epsilon}(x_0)}GdV_g-
  \dint_{\partial B_{L\epsilon}(x_0)}\frac{\partial G}{\partial r}\Delta_gGdS_g+
  \dint_{\partial B_{L\epsilon}}
   G\frac{\partial \Delta_gG}{\partial r}dS_g\\[\mv]
  &=&-\frac{1}{16\pi^2}+S_0-\frac{\log L\epsilon}{8\pi^2}+O(L\epsilon\log L\epsilon)
\end{array}$$
Now let us find a condition to have \;
$\int_M|\Delta_gf_\epsilon|^2dV_g=1$. By trivial calculations we can
see that it is equivalent to
\begin{equation*}
\frac{1}{32\pi^2c^2}\left(-\frac{5}{3}+2\log(1+\lambda
L^2)+32\pi^2S_0-4\log L\epsilon +O(\frac{1}{L^2})+O(L\epsilon\log
L\epsilon) \right)=1.
\end{equation*}
i.e.
\begin{equation*}
32\pi^2c^2=-\frac{5}{3}+2\log(1+\lambda L^2)+32\pi^2S_0-4\log
L\epsilon+O(\frac{1}{L^2}) +O(L\epsilon\log L\epsilon) .
\end{equation*}
Hence by (\ref{Lambda})  \;$\Lambda$\;take the following form
\begin{equation*}
\Lambda=\frac{10}{3}-64\pi^2S_0+O(\frac{1}{L^2})+O(L\epsilon\log
L\epsilon).
\end{equation*}
On the other hand it is easily seen that
\begin{equation*}
\int_{B_{L\epsilon}(x_0)}f_\epsilon dV_g=O(c(L\epsilon)^4);
\end{equation*}
and
\begin{equation*}
\int_{M\setminus B_{L\epsilon}(x_0)}f_\epsilon
dV_g=-\int_{B_{L\epsilon}}\frac{G}{c} =O(\frac{(L\epsilon)^4\log
L\epsilon}{c}).
\end{equation*}
hence
\begin{equation*}
\bar{f}_\epsilon=O(c(L\epsilon)^4).
\end{equation*}
Furthermore by trivial calculations one gets that in
\;$B_{L\epsilon}(x_0)$\;
$$\begin{array}{lll}
   (f_\epsilon-\bar{f}_\epsilon)^2&\geq& c^2+\frac{2}{64\pi^2}\left(\Lambda+B r^2
   -4\log(1+\lambda (\frac{r}{\epsilon})^2)+64\pi^2S_0+O(L\epsilon)
             +O(c^2(L\epsilon)^4)\right)\\[\mv]
          &=&c^2+\frac{5}{48\pi^2}-\frac{\log(1+\lambda (r/\epsilon)^2)}{8\pi^2}
          +O(\frac{1}{L^2})+O(L\epsilon\log L\epsilon)+O(c^2(L\epsilon)^4);
\end{array}$$
hence
$$\begin{array}{lll}
\dint_{B_{L\epsilon}(x_0)}e^{32\pi^2(f_\epsilon-\bar{f}_\epsilon)^2}dV_g&\geq&
 (1+O(L\epsilon)^2)
\dint_{B_{L\epsilon}(x_0)}e^{32\pi^2\left(c^2+\frac{5}{48\pi^2}-\frac{\log(1+\lambda
(r/\epsilon)^2}{8\pi^2}\right)+O(\frac{1}{L^2})+O(L\epsilon\log
L\epsilon)
+O(c^2(L\epsilon)^4)}dx\\[\mv]
&=&\epsilon^4e^{\frac{10}{3}+32\pi^2c^2+O(\frac{1}{L^2})+O(L\epsilon\log
L\epsilon) +O(c^2(L\epsilon)^4)}
\left(\pi^2\frac{L^6}{1+\lambda L^6}+O(L\epsilon)^2\right)\\[\mv]
&=&\epsilon^4e^{\frac{10}{3}+32\pi^2c^2}\pi^2
(1+O(\frac{1}{L^2})+O(L\epsilon\log L\epsilon)+O(L\epsilon)^2)\\[\mv]
&=&\frac{\pi^2}{6}e^{\frac{5}{3}+32\pi^2S_0}(1+O(L\epsilon\log
L\epsilon)+O(\frac{1}{L^2}) +O(c^2(L\epsilon)^4)).
\end{array}$$
on the other hand
$$\begin{array}{lll}
      \dint_{M\setminus B_{L\epsilon}(x_0)}e^{32\pi^2(f_\epsilon-\bar f_{\epsilon})^2}
      dV_g&\geq&\dint_{M\setminus B_{L\epsilon}(x_0)}(1+32\pi^2(f_\epsilon-\bar
          f_{\epsilon})^2)dV_g\\[\mv]
      &\geq&Vol(M\setminus B_{L\epsilon}(x_0))+
            \frac{\dint_{M\setminus B_{L\epsilon}(x_0)}32\pi^2G^2dV_g+O(c(L\epsilon)^4)}
            {c^2}\\[\mv]
     &=&Vol(M)+\frac{\dint_{M}32\pi^2G^2dV_g}{c^2}+O(L\epsilon)^4\log L\epsilon
\end{array}$$
Thus we arrive to
\begin{equation*}
\begin{split}
\int_Me^{32\pi^2(f_\epsilon-\bar{f}_\epsilon)^2}dV_g
   \geq Vol(M)+\frac{\pi^2}{6}e^{\frac{5}{3}+32\pi^2S_0}+\frac{\int_{M\setminus
   B_{L\epsilon}(x_0)}32\pi^2G^2dV_g}
    {c^2}\\+O(L\epsilon\log(L\epsilon))+O(\frac{1}{L^2})+O(c^2(L\epsilon)^4)
    \end{split};
    \end{equation*}
    and factorizing by \;$\frac{1}{c^2}$\;we get
    \begin{equation*}
  \begin{array}{lll}
\dint_Me^{32\pi^2(f_\epsilon-\bar{f}_\epsilon)^2}dV_g
   &\geq& Vol(M)+\frac{\pi^2}{6}e^{\frac{5}{3}+32\pi^2S_0}\\[\mv]
   &&+\frac{1}{c^2}\left(\dint_{M}32\pi^2G^2dV_g
    +O(c^2L\epsilon\log(L\epsilon))+O(\frac{c^2}{L^2})
    +O(c^4(L\epsilon)^4)\right).
\end{array}
\end{equation*}
On the other hand setting
\begin{equation*}
L=\log\frac{1}{\epsilon}
\end{equation*}
we get
\begin{equation*}
O(c^2L\epsilon\log(L\epsilon))+O(\frac{c^2}{L^2})+O(c^4(L\epsilon)^4)\rightarrow
0\;\; \text{as}\;\;\epsilon\rightarrow 0.
\end{equation*}
Hence the Proposition is proved.
\subsection{Proof of  Theorem\;$\ref{eq:theorem1}$}
This small subsection is concerned about the proof of Theorem\;$\ref{eq:theorem1}$.\\
First of all by corollary we have that
\begin{equation*}
\lim_{k\rightarrow +\infty}
\int_Me^{\alpha_ku_k^2}=Vol_g(M)+\tau^2\lim_{k\rightarrow+\infty}\frac{\lambda_k}{\beta_k^2}
\end{equation*}
with \;$\tau\neq 0$.\\
On the other hand from Proposition\;$\ref{eq:upperbound}$\; we get
\begin{equation*}
\tau^2\lim_{k\rightarrow+\infty}\frac{\lambda_k}{\beta_k^2} \leq
\frac{\pi^2}{6}e^{\frac{5}{3}+32\pi^2S_0}.
\end{equation*}
Hence we obtain
\begin{equation*}
\lim_{k\rightarrow +\infty} \int_Me^{\alpha_ku_k^2}\leq
Vol_g(M)+\frac{\pi^2}{6}e^{\frac{5}{3}+32\pi^2S_0}.
\end{equation*}
Thus using the relation
\begin{equation*}
\lim\limits_{k\rightarrow
+\infty}\int_{M}e^{\alpha_ku_k^2}dV_g=\sup_{u\in\mathcal {H}_1}
\int_{M}e^{32\pi^2u^2}dV_g.
\end{equation*}
we derive
\begin{equation*}
\sup_{u\in\mathcal {H}_1}\int_{M}e^{32\pi^2u^2}dV_g\leq
Vol_g(M)+\frac{\pi^2}{6}e^{\frac{5}{3} +32\pi^2S_0}.
\end{equation*}
On the other hand from Proposition\;$\ref{eq:test}$\; we have the
existence of a family of function \;$f_{\epsilon}$\;such that
\begin{equation*}
\int_M|\Delta_gf_\epsilon|^2dV_g=1;
\end{equation*}
and
\begin{equation*}
\limsup_{\epsilon\rightarrow
0}\int_Me^{32\pi^2(f_\epsilon-\bar{f}_\epsilon)^2}dV_g
   >Vol(M)+\frac{1}{6}e^{\frac{5}{3}+32\pi^2S_0}\pi^2.
\end{equation*}
Hence we reach a contradiction. So the proof of
Theorem\;$\ref{eq:theorem1}$\;is completed.
\end{proof}
\section{Proof of Theorem\;$\ref{eq:theorem2}$}
As already said in the Introduction, in this brief Section we will
explain how the proof
of Theorem\;$\ref{eq:theorem1}$\;remains valid for Theorem\;$\ref{eq:theorem2}$.\\
First of all we remark that all the analysis above have been
possible due to the following
facts\\
1)\\
\;$\int_M|\D_gu|^2dV_g$\;is an equivalent norm to the standard norm
of \;$H^2(M)$\;on\;
$\mathcal{H}_1$.\\
2)\\
The existence of the Green function for \;$\D_g^2$.\\
3)\\
The result of Fontana.\\
On the other hand we have a counterpart of 2) and 3). Moreover it is
easy to see that \;$\left<P^4_gu,u\right>$\;is also an equivalent
norm to the standard norm of \;$H^2(M)$\; on \;$\mathcal{H}_2$.
Notice that for a blowing-up sequence \;$u_k$\; we have that
\begin{equation}\label{eq:IMP}
\left<P^4_gu_k,u_k\right>=\int_{M}|\D_gu_k|^2dV_g+o_{k}(1);
\end{equation}
then it is easy to see that the same proof is valid up to the
subsection of test functions. Notice that \;$\eqref{eq:IMP}$\;holds
for the test functions \;$f_{\epsilon}$\;, then it is easy to see
that continuing the same proof we get Theorem\;$\ref{eq:theorem2}$.

\end{document}